\documentclass[a4paper,12pt]{amsart}
\usepackage[a4paper,hscale=0.7,vscale=0.75,centering]{geometry}
\usepackage{fullpage}
\usepackage{url}
\usepackage{nameref}
\usepackage{amsmath}
\usepackage{amssymb}

\newtheorem{theorem}{Theorem}[section]
\newtheorem{lemma}[theorem]{Lemma}
\newtheorem{corollary}[theorem]{Corollary}

\theoremstyle{definition}

\newtheorem{example}[theorem]{Example}

\newtheorem{assumption}{Assumption}

\newtheorem*{proofMainCoupling}{Proof of Theorem \ref{mainCouplingTheorem}}
\newtheorem*{proofDriftPerturbation}{Proof of Theorem \ref{theoremDriftPerturbation}}
\newtheorem*{proofGeneralMalliavin}{Proof of Theorem \ref{generalMalliavinResult}}
\newtheorem*{proofCorollaryMalliavin}{Proof of Corollary \ref{corollaryMalliavinBounds}}
\newtheorem*{proofTransportationInequalities}{Proof of Theorem \ref{transportationInequalitiesTheorem} and Theorem \ref{concentrationInequalitiesTheorem}}
\newtheorem*{proofCorollaryTransportation}{Proof of Corollary \ref{corollaryTransportation}}

\newenvironment{myassumption}[2][]
  {\renewcommand\theassumption{#2}\begin{assumption}[#1]}
  {\end{assumption}}

\theoremstyle{remark}
\newtheorem{remark}[theorem]{Remark}

\numberwithin{equation}{section}

\begin{document}

\title{Transportation inequalities for non-globally dissipative SDEs with jumps via Malliavin calculus and coupling}

\author{Mateusz B. Majka}
\address{Institute for Applied Mathematics, University of Bonn, Endenicher Allee 60, 53115 Bonn, Germany}

\email{majka@uni-bonn.de}



\date{}

\dedicatory{}

\begin{abstract}
By using the mirror coupling for solutions of SDEs driven by pure jump L\'{e}vy processes, we extend some transportation and concentration inequalities, which were previously known only in the case where the coefficients in the equation satisfy a~global dissipativity condition. Furthermore, by using the mirror coupling for the jump part and the coupling by reflection for the Brownian part, we extend analogous results for jump diffusions. To this end, we improve some previous results concerning such couplings and show how to combine the jump and the Brownian case. As a crucial step in our proof, we develop a novel method of bounding Malliavin derivatives of solutions of SDEs with both jump and Gaussian noise, which involves the coupling technique and which might be of independent interest. The bounds we obtain are new even in the case of diffusions without jumps.
\end{abstract}

\maketitle

\section{Introduction}

We consider stochastic differential equations in $\mathbb{R}^d$ of the form
\begin{equation}\label{SDE1}
 dX_t = b(X_t)dt + \sigma(X_t) dW_t + \int_{U} g(X_{t-},u) \widetilde{N}(dt,du) \,,
\end{equation}
where $(W_t)_{t \geq 0}$ is an $m$-dimensional Brownian motion and $\widetilde{N}(dt,du) = N(dt,du) - dt \,\nu(du)$ is a compensated Poisson random measure on $\mathbb{R}_{+} \times U$, where $(U,\mathcal{U}, \nu)$ is a $\sigma$-finite measure space. Moreover, the coefficients $b: \mathbb{R}^d \to \mathbb{R}^d$, $\sigma: \mathbb{R}^d \to \mathbb{R}^{d \times m}$ and $g: \mathbb{R}^d \times U \to \mathbb{R}^d$ are such that for any $x \in \mathbb{R}^d$ we have
\begin{equation*}
 \int_{U} |g(x,u)|^2 \nu(du) < \infty
\end{equation*}
and there exists a continuous function $\kappa : \mathbb{R}_{+} \to \mathbb{R}$ such that for all $x$, $y \in \mathbb{R}^d$ we have
\begin{equation}\label{onesidedLipschitz}
 \langle b(x) - b(y) , x - y \rangle + \frac{1}{2}\int_U |g(x,u) - g(y,u)|^2 \nu(du) + \| \sigma(x) - \sigma(y) \|_{HS}^2 \leq -\kappa(|x-y|)|x-y|^2 \,,
\end{equation}
where $\| \sigma \|_{HS} = \sqrt{\operatorname{tr} \sigma \sigma^T}$ is the Hilbert-Schmidt norm. Note that $\kappa$ is allowed to take negative values.

If the condition (\ref{onesidedLipschitz}) holds with a constant function $\kappa \equiv K$ for some $K \in \mathbb{R}$, we call (\ref{onesidedLipschitz}) a one-sided Lipschitz condition. If $K > 0$, we call it a (global) dissipativity condition. If a one-sided Lipschitz condition is satisfied and we additionally assume that the drift $b$ is continuous and that $\sigma$~and $g$ satisfy a linear growth condition, we can prove that (\ref{SDE1}) has a unique non-explosive strong solution, even if the one-sided Lipschitz condition is satisfied only locally (see e.g. Theorem 2 in \cite{gyongykrylov}).

 For $p \geq 1$, the $L^p$-\emph{Wasserstein distance} (or the $L^p$-transportation cost) between two probability measures $\mu_1$, $\mu_2$ on a metric space $(E, \rho)$ is defined by
 \begin{equation*}
  W_{p,\rho}(\mu_1, \mu_2) := \inf_{\pi \in \Pi(\mu_1, \mu_2)} \left( \int \int \rho(x,y)^p \pi(dx \, dy) \right)^{1/p} \,,
 \end{equation*}
 where $\Pi(\mu_1, \mu_2)$ is the family of all couplings of $\mu_1$ and $\mu_2$, i.e., $\pi \in \Pi(\mu_1, \mu_2)$ if and only if $\pi$ is a measure on $E \times E$ with marginals $\mu_1$ and $\mu_2$. If the metric space $(E,\rho)$ is chosen to be $\mathbb{R}^d$ with the Euclidean metric $\rho(x,y) = |x-y|$, then we denote $W_{p,\rho}$ just by $W_p$.
 
 If the equation (\ref{SDE1}) is globally dissipative with some constant $K > 0$, then it is well known that the solution $(X_t)_{t \geq 0}$ to (\ref{SDE1}) has an invariant measure and that the transition semigroup $(p_t)_{t \geq 0}$ associated with $(X_t)_{t \geq 0}$ is exponentially contractive with respect to $W_p$ for any $p \in [1,2]$, i.e.,
 \begin{equation*}
  W_p(\mu_1 p_t, \mu_2 p_t) \leq e^{-Kt} W_p(\mu_1, \mu_2)
 \end{equation*}
for any probability measures $\mu_1$ and $\mu_2$ on $\mathbb{R}^d$ and any $t > 0$ (see e.g. the proof of Theorem 2.2 in \cite{yma}). However, we will show that for $p = 1$ a related result still holds (under some additional assumptions, see Corollary \ref{corollaryExponentialContractivity}) if we replace the global dissipativity condition with the following one.
\begin{myassumption}{D1}\label{AssumptionD1} (Dissipativity at infinity)
\begin{equation*}
  \limsup_{r \to \infty} \kappa(r) > 0 \,.
\end{equation*}
\end{myassumption}
In other words, Assumption \ref{AssumptionD1} states that there exist constants $R > 0$ and $K > 0$ such that for all $x$, $y \in \mathbb{R}^d$ with $|x-y| > R$ we have
\begin{equation*}
 \langle b(x) - b(y) , x - y \rangle + \frac{1}{2}\int_U |g(x,u) - g(y,u)|^2 \nu(du) + \| \sigma(x) - \sigma(y) \|_{HS}^2 \leq -K|x-y|^2 \,,
\end{equation*}
which justifies calling it a dissipativity at infinity condition. Moreover, in some cases we will also need another condition on the function $\kappa$, namely
\begin{myassumption}{D2}\label{AssumptionD2} (Regularity of the drift at zero)
\begin{equation*}
  \lim_{r \to 0} r \kappa(r) = 0 \,.
\end{equation*}
\end{myassumption}
This is obviously satisfied if, e.g., there is a constant $L > 0$ such that we have $\kappa(r) \geq - L$ for all $r \geq 0$ (which is the case whenever the coefficients in (\ref{SDE1}) satisfy a one-sided Lipschitz condition) and if $b$ is continuous. Such an assumption is quite natural in order to ensure existence of a solution to (\ref{SDE1}).

 For probability measures $\mu_1$ and $\mu_2$ on $(E,\rho)$, we define the \emph{relative entropy} (Kullback-Leibler information) of $\mu_1$ with respect to $\mu_2$ by
\begin{equation*}
 H(\mu_1 | \mu_2) := \begin{cases}
                  \int \log \frac{d \mu_1}{d \mu_2} d \mu_1 & \text{if } \mu_1 \ll \mu_2 \,, \\
                  +\infty & \text{otherwise} \,.
                 \end{cases}
\end{equation*}

We say that a probability measure $\mu$ satisfies an $L^p$-\emph{transportation cost-information inequality} on $(E,\rho)$ if there is a constant $C > 0$ such that for any probability measure $\eta$~we have
\begin{equation*}
  W_{p,\rho}(\eta, \mu) \leq \sqrt{2C H(\eta | \mu)} \,.
\end{equation*}
Then we write $\mu \in T_p(C)$.

The most important cases are $p = 1$ and $p = 2$. Since $W_{1,\rho} \leq W_{2,\rho}$, we see that the $L^2$-transportation inequality (the $T_2$ inequality, also known as the Talagrand inequality) implies $T_1$, and it is well known that in fact $T_2$ is much stronger. The $T_2$ inequality has some interesting connections with other well-known functional inequalities. Due to Otto and Villani \cite{ottovillani}, we know that the log-Sobolev inequality implies $T_2$, whereas $T_2$ implies the Poincar\'{e} inequality. On the other hand, the $T_1$ inequality is related to the phenomenon of measure concentration. Indeed, consider a generalization of $T_1$ known as the $\alpha$-$W_1H$ inequality. Namely, let $\alpha$ be a non-decreasing, left continuous function on $\mathbb{R}_{+}$ with $\alpha(0)=0$. We say that a probability measure $\mu$ satisfies a $W_1H$-inequality with deviation function $\alpha$ (or simply $\alpha$-$W_1 H$ inequality) if for any probability measure $\eta$ we have
 \begin{equation}\label{alphaW1H}
  \alpha(W_{1,\rho}(\eta, \mu)) \leq H(\eta | \mu) \,.
 \end{equation}
 
 We have the following result which is due to Gozlan and L\'{e}onard (see Theorem 2 in \cite{gozlan} for the original result, cf. also Lemma 2.1 in \cite{lwu}). It is a generalization of a result by Bobkov and G\"{o}tze (Theorem 3.1 in \cite{bobkov}), which held only for the quadratic deviation function.

 Fix a probability measure $\mu$ on $(E,\rho)$ and a convex deviation function $\alpha$. Then the following properties are equivalent:
\begin{enumerate}
 \item the $\alpha$-$W_1 H$ inequality for the measure $\mu$ holds, i.e., for any probability measure $\eta$ on $(E,\rho)$ we have
 \begin{equation*}
  \alpha(W_{1,\rho}(\eta, \mu)) \leq H(\eta | \mu) \,,
 \end{equation*}
\item for every $f : E \to \mathbb{R}$ bounded and Lipschitz with $\| f \|_{\operatorname{Lip}} \leq 1$ we have
\begin{equation}\label{gozlanCharacterization}
 \int e^{\lambda (f - \mu(f))} d\mu \leq e^{\alpha^*(\lambda)} \text{ for any } \lambda > 0 \,,
\end{equation}
where $\alpha^*(\lambda) := \sup_{r \geq 0}(r\lambda - \alpha (r))$ is the convex conjugate of $\alpha$,
\item if $(\xi_k)_{k \geq 1}$ is a sequence of i.i.d random variables with common law $\mu$, then for every $f : E \to \mathbb{R}$ bounded and Lipschitz with $\| f \|_{\operatorname{Lip}} \leq 1$ we have
\begin{equation}\label{probabilisticCharacterization}
 \mathbb{P} \left( \frac{1}{n} \sum_{k=1}^{n} f(\xi_k) - \mu(f) > r \right) \leq e^{-n \alpha(r)} \text{ for any } r > 0, n \geq 1 \,.
\end{equation}
\end{enumerate}
 
 This gives an intuitive interpretation of $\alpha$-$W_1 H$ in terms of a concentration of measure property (\ref{probabilisticCharacterization}), while the second characterization (\ref{gozlanCharacterization}) is very useful for proving such inequalities, as we shall see in the sequel. For a general survey of transportation inequalities the reader might consult \cite{gozlansurvey} or Chapter 22 of \cite{villani}.

As an example of a simple equation of the type (\ref{SDE1}) consider
\begin{equation*}
 dX_t = b(X_t)dt + \sqrt{2} dW_t
\end{equation*}
with a $d$-dimensional Brownian motion $(W_t)_{t \geq 0}$. If the global dissipativity assumption is satisfied, then $(X_t)_{t \geq 0}$ has an invariant measure $\mu$ and by a result of Bakry and \'{E}mery \cite{bakryemery}, $\mu$ satisfies the log-Sobolev inequality and thus (by Otto and Villani \cite{ottovillani}) also the Talagrand inequality.
More generally, for equations of the form
\begin{equation}\label{GaussianMultiplicativeSDE}
 dX_t = b(X_t)dt + \sigma(X_t)dW_t \,,
\end{equation}
also under the global dissipativity assumption, Djellout, Guillin and Wu in \cite{djellout} showed that $T_2$ holds for the invariant measure, as well as on the path space. As far as we are aware, there are currently no results in the literature concerning transportation inequalities for equations like (\ref{GaussianMultiplicativeSDE}) without assuming global dissipativity. Hence, even though in the present paper we focus on SDEs with jumps, our results may be also new in the purely Gaussian case.

 For equations of the form
 \begin{equation}\label{PoissonJumpSDE}
  dX_t = b(X_t)dt + \int_{U} g(X_{t-},u) \widetilde{N}(dt,du) \,,
 \end{equation}
the Poincar\'{e} inequality does not always hold (see Example 1.1 in \cite{lwu}) and thus in general we cannot have $T_2$. However, under the global dissipativity assumption, Wu in \cite{lwu} showed some $\alpha$-$W_1H$ inequalities.

 Suppose there is a real measurable function $g_{\infty}$ on $U$ such that $|g(x,u)| \leq g_{\infty}(u)$ for every $x \in \mathbb{R}^d$ and $u \in U$. We make the following assumption.
 \begin{myassumption}{E}\label{AssumptionE}
  (Exponential integrability of the intensity measure)\\
  There exists a constant $\lambda > 0$ such that
  \begin{equation*}
      \beta(\lambda) := \int_{U} (e^{\lambda g_{\infty}(u) } - \lambda g_{\infty}(u) - 1) \nu(du) < \infty \,,
  \end{equation*}
  where $\nu$ is the intensity measure associated with $N$.
 \end{myassumption}

 \begin{remark}
  Assumption \ref{AssumptionE} is quite restrictive. In particular, let us consider the case where $U \subset \mathbb{R}^d$ and $g(x,u) = \widetilde{g}(x) u$ for some $\mathbb{R}^{d \times d}$-valued function $\widetilde{g}$ and hence the equation (\ref{PoissonJumpSDE}) is driven by a $d$-dimensional L\'{e}vy process $(L_t)_{t \geq 0}$ (i.e., we have $dX_t = b(X_t)dt + \widetilde{g}(X_{t-}) dL_t$). Then Assumption \ref{AssumptionE} implies finiteness of an exponential moment of $(L_t)_{t \geq 0}$ (cf. Theorem 25.3 and Corollary 25.8 in \cite{sato}). However, there are examples of equations of such type for which the $\alpha$-$W_1H$ inequality implies Assumption \ref{AssumptionE}, and hence in general we cannot prove such inequalities without it (see Remark 2.5 in \cite{lwu}). Nevertheless, without this assumption it is still possible to obtain some concentration inequalities (see Remark 5.2 in \cite{lwu} or Theorem \ref{concentrationInequalitiesTheorem} below).
 \end{remark}

 Fix $T > 0$ and define a deviation function
\begin{equation*}
 \alpha_{T}(r) := \sup_{\lambda \geq 0} \left\{ r \lambda - \int_0^T \beta(e^{-Kt}\lambda)dt \right\} \,,
\end{equation*}
where the constants $\lambda > 0$ and $K > 0$ are such that Assumption \ref{AssumptionE} is satisfied with $\lambda$ and that (\ref{PoissonJumpSDE}) is globally dissipative with the dissipativity constant $K$. Then for any $T > 0$ and any $x \in \mathbb{R}^d$, by Theorem 2.2 in \cite{lwu} we have the $W_1H$ transportation inequality with deviation function $\alpha_T$ for the measure $\delta_x p_T$, which is the law of the random variable $X_T(x)$, where $(X_t(x))_{t \geq 0}$ is a solution to (\ref{PoissonJumpSDE}) starting from $x \in \mathbb{R}^d$, i.e., we have
 \begin{equation*}
  \alpha_T(W_1(\eta, \delta_x p_T)) \leq H(\eta | \delta_x p_T)
 \end{equation*}
for any probability measure $\eta$ on $\mathbb{R}^d$, where $W_1 = W_{1,\rho}$ with $\rho$ being the Euclidean metric on $\mathbb{R}^d$. Analogous results have been proved by a very similar approach in \cite{yma} for equations of the form (\ref{SDE1}), i.e., including also the Gaussian noise. 

In the sequel we will explain how to modify the proofs in \cite{lwu} and \cite{yma} to replace the global dissipativity assumption with our Assumption \ref{AssumptionD1}. We will show that we can obtain $\alpha$-$W_1H$ inequalities by using couplings to control perturbations of solutions to (\ref{SDE1}), see Theorem \ref{transportationInequalitiesTheorem}. We will also prove that the construction of the required couplings is possible for a~certain class of equations satisfying Assumption \ref{AssumptionD1} (Theorems \ref{mainCouplingTheorem} and \ref{theoremDriftPerturbation}). All these results together will imply our extension of the main theorems from \cite{lwu} and \cite{yma}, which is stated as Corollary \ref{corollaryTransportation}.

 The method of the proof is based on the Malliavin calculus. On any filtered probability space $(\Omega, \mathcal{F}, (\mathcal{F}_t)_{t \geq 0}, \mathbb{P})$ equipped with an $m$-dimensional Brownian motion $(W_t)_{t \geq 0}$ and a Poisson random measure $N$ on $\mathbb{R}_{+} \times U$, we can define the Malliavin derivatives for a~certain class of measurable functionals $F$ with respect to the process $(W_t)_{t \geq 0}$ (the classic Malliavin differential operator $\nabla$), as well as a Malliavin derivative of $F$ with respect to $N$ (the difference operator $D$). Namely, if we consider the family $\mathcal{S}$ of smooth functionals of $(W_t)_{t \geq 0}$ of the form
\begin{equation*}
 F = f(W(h_1), \ldots , W(h_n)) \text{ for } n \geq 1 \,,
\end{equation*}
where $W(h) = \int_0^T h(s) dW_s$ for $h \in H = L^2([0,T] ; \mathbb{R}^m)$ and $f \in C^{\infty}(\mathbb{R}^n)$, we can define the Malliavin derivative with respect to $(W_t)_{t \geq 0}$ as the unique element $\nabla F$ in $L^2(\Omega ; H) \simeq L^2(\Omega \times [0,T] ; \mathbb{R}^m)$ such that for any $h \in H$ we have
\begin{equation*}
\langle \nabla F , h \rangle_{L^2([0,T];\mathbb{R}^m)} = \lim_{\varepsilon \to 0} \frac{1}{\varepsilon} \left( F(W_{\cdot} + \int_0^{\cdot} h_s ds) - F(W_{\cdot}) \right) \,,
\end{equation*}
where the convergence holds in $L^2(\Omega)$ (see e.g. Definition A.10 in \cite{dinunno}). Then the definition can be extended to all random variables $F$ in the space $\mathbb{D}^{1,2}$ which is the completion of $\mathcal{S}$ in $L^2(\Omega)$ with respect to the norm
\begin{equation*}
\| F \|_{\mathbb{D}^{1,2}}^2 := \| F \|_{L^2(\Omega)}^2 + \| \nabla F \|_{L^2(\Omega ; H)}^2 \,.
\end{equation*}
For a brief introduction to the Malliavin calculus with respect to Brownian motion see Appendix A in \cite{dinunno} or Chapter VIII in \cite{bass} and for a comprehensive treatment the monograph \cite{nualart}.
On the other hand, the definition of the Malliavin derivative with respect to $N$ that we need is much less technical, since it is just a difference operator. Namely, if our Poisson random measure $N$ on $\mathbb{R}_{+} \times U$ has the form
\begin{equation*}
N = \sum_{j=1}^{\infty} \delta_{(\tau_j, \xi_j)}
\end{equation*}
with $\mathbb{R}_{+}$-valued random variables $\tau_j$ and $U$-valued $\xi_j$, then for any measurable functional $f$ of $N$ and for any $(t,u) \in \mathbb{R}_{+} \times U$ we put
\begin{equation}\label{defDifferenceOperator}
D_{t,u} f(N) := f(N + \delta_{(t,u)}) - f(N) \,.
\end{equation}
There is also an alternative approach to the Malliavin calculus for jump processes, where the Malliavin derivative is defined as an actual differential operator, which was in fact the original approach and which traces back to Bismut \cite{bismut}, see also \cite{basscranston} and \cite{bichteler}. However, for our purposes we prefer the definition (\ref{defDifferenceOperator}), which was introduced by Picard in \cite{picardfrench} and \cite{picard}, and which is suitable for proving the Clark-Ocone formula. Namely, we will need to use the result stating that for any $F$ being a functional of $(W_t)_{t \geq 0}$ and $N$ such that
 \begin{equation}\label{integrabilityConditionForClarkOcone}
  \mathbb{E} \int_0^T |\nabla_t F|^2 dt + \mathbb{E} \int_0^T \int_U |D_{t,u}F|^2 \nu(du) dt < \infty \,,
 \end{equation}
we have
\begin{equation*}
 F = \mathbb{E}F + \int_0^T \mathbb{E}[\nabla_t F | \mathcal{F}_t] dW_t + \int_0^T \int_U \mathbb{E}[D_{t,u} F | \mathcal{F}_t] \widetilde{N}(dt,du) \,.
\end{equation*}
 It is proved in \cite{lokka} that the definition (\ref{defDifferenceOperator}) is actually equivalent to the definition of the Malliavin derivative for jump processes via the chaos expansion and this approach is used to obtain the Clark-Ocone formula for the pure jump case. For the jump diffusion case, see Theorem 12.20 in \cite{dinunno}. For more general recent extensions of this result, see \cite{lastpenrose}. Once we apply the Clark-Ocone formula to the solution of (\ref{SDE1}), we can obtain some information on its behaviour by controlling its Malliavin derivatives. Therefore one of the crucial components of the proof of our results in this paper is Theorem \ref{generalMalliavinResult}, presenting a novel method of bounding such derivatives, which, contrary to the method used in Lemma 3.4 in \cite{yma}, works also without the global dissipativity assumption and without any explicit regularity conditions on the coefficients of (\ref{SDE1}), except some sufficient ones to guarantee Malliavin differentiability of the solution (it is enough if the coefficients are Lipschitz, see e.g. Theorem 17.4 in \cite{dinunno}).

The last notion that we need to introduce before we will be able to formulate our main results is that of a coupling. For an $\mathbb{R}^d$-valued Markov process $(X_t)_{t \geq 0}$ with transition kernels $(p_t(x,\cdot))_{t \geq 0, x \in \mathbb{R}^d}$ we say that an $\mathbb{R}^{2d}$-valued process $(X_t',X_t'')_{t \geq 0}$ is a \emph{coupling} of two copies of the Markov process $(X_t)_{t \geq 0}$ if both $(X_t')_{t \geq 0}$ and $(X_t'')_{t \geq 0}$ are Markov processes with transition kernels $p_t$ but possibly with different initial distributions. The construction of appropriate couplings of solutions to equations like (\ref{SDE1}) plays the key role in the proofs of Theorems \ref{mainCouplingTheorem} and \ref{theoremDriftPerturbation}. For more information about couplings, see e.g. \cite{lindvall}, \cite{eberle}, \cite{majka} and the references therein.

The only papers that we are aware of which deal with transportation inequalities directly in the context of SDEs with jumps are \cite{lwu}, \cite{maprivault}, \cite{yma} and \cite{shao}. The latter two actually extend the method developed by Wu in \cite{lwu}, but in both these papers a kind of global dissipativity assumption is required (see Remark \ref{remarkShao} for a discussion about \cite{shao}). In the present paper we explain how to drop this assumption (by imposing some additional conditions) and further extend the method of Wu. Since our extension lies at the very core of the method, it allows us to improve on essentially all the main results and corollaries obtained in \cite{lwu} and \cite{yma} (and it might be also applicable to the results in \cite{shao}, cf. once again Remark \ref{remarkShao}), replacing the global dissipativity assumption with a~weaker condition.

On the other hand, in \cite{maprivault} some convex concentration inequalities of the type (\ref{concentrationIneq1}) have been shown for a certain class of additive functionals $S_T = \int_0^T g(X_t) dt$ of solutions $(X_t)_{t \geq 0}$ to equations like (\ref{SDE1}). These are later used to obtain some $\alpha$-$W_1I$ inequalities, which are analogous to $\alpha$-$W_1H$ inequalities (\ref{alphaW1H}) but with the Kullback-Leibler information $H$~replaced with the Fisher-Donsker-Varadhan information, see e.g. \cite{guillinleonard} for more details. The proof in \cite{maprivault}, similarly to \cite{lwu}, is based on the forward-backward martingale method from \cite{klein}, but unlike \cite{lwu} it does not use the Malliavin calculus. In the framework of Wu from \cite{lwu} that we use here, it is possible to obtain related $\alpha$-$W_1J$ inequalities with $J$ being the modified Donsker-Varadhan information. Once we have transportation inequalities like the ones in our Theorem \ref{transportationInequalitiesTheorem}, we can use the methods from Corollary 2.15 in \cite{lwu} and Corollary 2.7 in \cite{yma}. This is, however, beyond the scope of the present paper and in the sequel we focus on extending the main results from \cite{lwu} and \cite{yma}.

 \section{Main results}

We start with a general theorem, which shows that a key tool to obtain transportation inequalities for a solution $(X_t)_{t \geq 0}$ to 
\begin{equation}\label{generalSDEsect2}
   dX_t = b(X_t)dt + \sigma(X_t) dW_t + \int_{U} g(X_{t-},u) \widetilde{N}(dt,du)
  \end{equation}
is to be able to control perturbations of $(X_t)_{t \geq 0}$ via a coupling, with respect to changes in initial conditions (see (\ref{initialChangeBound}) below) as well as changes of the drift (see (\ref{driftChangeBound})). In the next two theorems we assume that the coefficients in (\ref{generalSDEsect2}) satisfy some sufficient conditions for existence of a solution and its Malliavin differentiability (e.g. they are Lipschitz, cf. Theorem 17.4 in \cite{dinunno}). From now on, $(\mathcal{F}_t)_{t \geq 0}$ will always denote the filtration generated by all the sources of noise in the equations that we consider, while $(p_t)_{t \geq 0}$ will be the transition semigroup associated with the solution to the equation. Moreover, for a process $(h_t)_{t \geq 0}$ adapted to $(\mathcal{F}_t)_{t \geq 0}$, we will denote by $(\widetilde{X}_t)_{t \geq 0}$ a solution to
  \begin{equation}\label{perturbedSDEsect2}
   d\widetilde{X}_t = b(\widetilde{X}_t)dt + \sigma(\widetilde{X}_t) h_t dt + \sigma(\widetilde{X}_t) dW_t + \int_{U} g(\widetilde{X}_{t-},u) \widetilde{N}(dt,du) \,.
  \end{equation}
 Then we have the following result.

  \begin{theorem}\label{transportationInequalitiesTheorem}
 Assume there exists a constant $\sigma_{\infty}$ such that for any $x \in \mathbb{R}^d$ we have $\| \sigma(x) \| \leq \sigma_{\infty}$, where $\| \cdot \|$ is the operator norm, and there exists a measurable function $g_{\infty}: U \to \mathbb{R}$ such that $|g(x,u)| \leq g_{\infty}(u)$ for any $x \in \mathbb{R}^d$ and $u \in U$. Assume further that there exists some $\lambda > 0$ such that Assumption \ref{AssumptionE} is satisfied. Moreover, suppose that there exists a coupling $(X_t,Y_t)_{t \geq 0}$ of solutions to (\ref{generalSDEsect2}) and a function $c_1 :\mathbb{R}_{+} \to \mathbb{R}_{+}$ such that for any $0 \leq s \leq t$ we have
 \begin{equation}\label{initialChangeBound}
  \mathbb{E}[|X_t - Y_t| / \mathcal{F}_s ]\leq c_1(t-s)|X_s-Y_s| \,.
 \end{equation}
  Furthermore, assume that there exists a coupling $(X_t,Y'_t)_{t \geq 0}$ of solutions to (\ref{generalSDEsect2}) and functions $c_2$, $c_3 :\mathbb{R}_{+} \to \mathbb{R}_{+}$ such that for any $0 \leq s \leq t$ we have
  \begin{equation}\label{driftChangeBound}
   \mathbb{E}[|\widetilde{X}_t - Y'_t| / \mathcal{F}_s ] \leq c_2(t-s) \mathbb{E} \int_s^t c_3(r) |\sigma(\widetilde{X}_r)h_r| dr \,.
  \end{equation}
 Then the following assertions hold.
 \begin{enumerate}
  \item For any $T > 0$ and for any $x \in \mathbb{R}^d$ the measure $\delta_x p_T$ satisfies
  \begin{equation}\label{transIneq1}
   \alpha_T(W_1(\eta, \delta_x p_T)) \leq H(\eta | \delta_x p_T)
  \end{equation}
for any probability measure $\eta$ on $\mathbb{R}^d$. Here $W_1 = W_{1,\rho}$ with $\rho$ being the Euclidean metric on $\mathbb{R}^d$ and
\begin{equation*}
 \alpha_{T}(r) := \sup_{\lambda \geq 0} \left\{ r \lambda - \int_0^T \beta(c_1(T-t) \lambda)dt - \frac{\sigma_{\infty}^2 c^2_2(T) \lambda^2}{2} \int_0^T c_3^2(t)dt \right\} \,.
\end{equation*}
\item For any $T > 0$ and for any $x \in \mathbb{R}^d$ the law $\mathbb{P}_{x,[0,T]}$ of $(X_t(x))_{t \in [0,T]}$ as a measure on the space $\mathbb{D}([0,T]; \mathbb{R}^d)$ of c\'{a}dl\'{a}g $\mathbb{R}^d$-valued functions on $[0,T]$ satisfies
\begin{equation}\label{transIneq2}
  \alpha^P_T(W_{1,d_{L^1}}(Q, \mathbb{P}_{x,[0,T]})) \leq H(Q | \mathbb{P}_{x,[0,T]})
\end{equation}
for any probability measure $Q$ on $\mathbb{D}([0,T]; \mathbb{R}^d)$. Here we take $d_{L^1}(\gamma_1, \gamma_2) := \int_0^T |\gamma_1(t) - \gamma_2(t)| dt$ as the $L^1$ metric on the path space and
\begin{equation*}
 \alpha^P_T(r) := \sup_{\lambda \geq 0} \left\{ r \lambda - \int_0^T \beta\left(\lambda \int_t^T c_1(s-t)ds  \right)dt - \frac{\sigma_{\infty}^2 \lambda^2}{2} \int_0^T c_3^2(t) \left( \int_t^T c_2(r)dr \right)^2 dt \right\} \,.
\end{equation*}
 \end{enumerate}
 \end{theorem}
 
 Note that in (\ref{driftChangeBound}) the process $(Y'_t)_{t \geq 0}$ is coupled with $(X_t)_{t \geq 0}$, but the estimated distance is between $(Y'_t)_{t \geq 0}$ and $(\widetilde{X}_t)_{t \geq 0}$ given by (\ref{perturbedSDEsect2}). In other words, we need to consider a process $(Y'_t)_{t \geq 0}$ whose law is determined by the dynamics (\ref{generalSDEsect2}), but it behaves in such a way that we can control its distance to a process following the modified dynamics given by (\ref{perturbedSDEsect2}). An example of such a construction can be found in the proof of Theorem \ref{theoremDriftPerturbation}.
 
 Even without Assumption \ref{AssumptionE}, it is still possible to recover some concentration inequalities.
 \begin{theorem}\label{concentrationInequalitiesTheorem}
  Assume that all the assumptions of Theorem \ref{transportationInequalitiesTheorem} are satisfied except for Assumption \ref{AssumptionE}. Instead, suppose that $g_{\infty}(u)$ is just square integrable with respect to $\nu$. Fix any $T > 0$ and any $x \in \mathbb{R}^d$. Then for any $C^2$ convex function $\phi$~such that $\phi'$ is also convex and for any Lipschitz function $f: \mathbb{R}^d \to \mathbb{R}$, we have
  \begin{equation}\label{concentrationIneq1}
  \begin{split}
   \mathbb{E}&\phi\Big(f(X_T(x)) - p_Tf(x)\Big) \\
   &\leq \mathbb{E}\phi\left(\| f \|_{\operatorname{Lip}} \left(\int_0^T \int_U c_1(T-t)g_{\infty}(u) \widetilde{N}(dt,du) + c_2(T) \int_0^T c_3(t) j(t)dW_t \right)\right) \,,
   \end{split}
  \end{equation}
where $j$ is any deterministic $\mathbb{R}^m$-valued function such that for all $t > 0$ we have $|j(t)| = \sigma_{\infty}$. Moreover, for any Lipschitz function $F: \mathbb{D}([0,T];\mathbb{R}^d) \to \mathbb{R}$ we have
\begin{equation}\label{concentrationIneq2}
\begin{split}
 \mathbb{E}&\phi\Big(F(X_{[0,T]}(x)) - \mathbb{E}F(X_{[0,T]}(x))\Big) \\
 &\leq \mathbb{E}\phi\Bigg(\| F \|_{\operatorname{Lip}} \Bigg(\int_0^T \int_U \left(\int_t^T c_1(r-t)dr\right) g_{\infty}(u) \widetilde{N}(dt,du) \\
 &+ \int_0^T c_3(t) \left(\int_t^T c_2(r)dr\right) j(t) dW_t \Bigg)\Bigg) \,.
 \end{split}
\end{equation}

 \end{theorem}

 The crucial step in proving the above theorems is to find appropriate bounds on Malliavin derivatives of the solution to (\ref{generalSDEsect2}). We will show that we can obtain such bounds on $D$ and $\nabla$ using conditions (\ref{initialChangeBound}) and (\ref{driftChangeBound}), respectively (see Section \ref{sectionMalliavin} for details).
 
 Now we present another result, which will consequently lead us to some examples of equations for which the inequalities (\ref{initialChangeBound}) and (\ref{driftChangeBound}) actually hold. First, however, we need to formulate some additional assumptions. We will need a pure jump L\'{e}vy process $(L_t)_{t \geq 0}$ with a L\'{e}vy measure $\nu^L$ satisfying the following set of conditions.
 \begin{myassumption}{L1}\label{Assumption1}
 (Rotational invariance of the L\'{e}vy measure)
 $\nu^L$ is rotationally invariant, i.e.,
 \begin{equation*}
 \nu^L(AB) = \nu^L(B)
\end{equation*}
for every Borel set $B \in \mathcal{B}(\mathbb{R}^d)$ and every $d \times d$ orthogonal matrix $A$.
\end{myassumption}

\begin{myassumption}{L2}\label{Assumption2}
(Absolute continuity of the L\'{e}vy measure) 
$\nu^L$ is absolutely continuous with respect to the Lebesgue measure on $\mathbb{R}^d$ with a density $q$ that is continuous almost everywhere on $\mathbb{R}^d$.
\end{myassumption}

 Under Assumptions \ref{Assumption1}-\ref{Assumption2} it has been proved in \cite{majka} (see Theorem 1.1 therein) that there exists a coupling $(X_t,Y_t)_{t \geq 0}$ of solutions to
 \begin{equation*}
  dX_t = b(X_t)dt + dL_t \,,
 \end{equation*}
 defined as a unique strong solution to the $2d$-dimensional SDE given in the sequel by (\ref{eqX}) and (\ref{optimalY}). Moreover, consider two additional conditions on the jump density $q$.
 
 \begin{myassumption}{L3}\label{Assumption3}
 (Positive mass of the overlap of the jump density and its translation)
 There exist constants $m$, $\delta > 0$ such that $\delta < 2m$ and
 \begin{equation}\label{overlapCondition}
  \inf_{x \in \mathbb{R}^d : 0 < |x| \leq \delta} \int_{\{ |v| \leq m \} \cap \{ |v + x| \leq m \}} q(v) \wedge q(v + x) dv > 0 \,.
 \end{equation}
\end{myassumption}
 
 \begin{myassumption}{L4}\label{Assumption4}
(Positive mass in a neighbourhood of zero)
 There exists a constant $\varepsilon > 0$ such that $\varepsilon \leq \delta$ (with $\delta$ defined via (\ref{overlapCondition}) above) and
 \begin{equation*}
  \int_{\{ |v| \leq \varepsilon / 2 \}} q(v) dv > 0 \,.
 \end{equation*}
\end{myassumption}
 
 Suppose now that all the Assumptions \ref{Assumption1}-\ref{Assumption4} are satisfied. Let us define a continuous function $\kappa : \mathbb{R}_{+} \to \mathbb{R}$ so that for any $x$, $y \in \mathbb{R}^d$ the condition $\langle b(x) - b(y) , x -y \rangle \leq -\kappa(|x-y|)|x-y|^2$ is satisfied and suppose that Assumption \ref{AssumptionD1} holds. Then we get that, by the inequality (1.8) in Theorem 1.1 in \cite{majka}, there exist explicitly given $L$, $\theta > 0$ and a~function $f : \mathbb{R}_{+} \to \mathbb{R}_{+}$ such that
\begin{equation}\label{majkaThm12}
 \mathbb{E}|X_t(x) - Y_t(y)| \leq L e^{-\theta t}f(|x-y|) \,.
\end{equation}
However, the function $f$~used in \cite{majka} is discontinuous. It is actually of the form
\begin{equation}\label{formulaF}
 f = a\mathbf{1}_{(0, \infty)} + f_1
\end{equation}
with $a > 0$ and $f_1$ being a continuous, concave function, extended in an affine way from some point $R_1 > 0$ (and thus we have $a_1x \leq f_1(x) \leq a_2x$ for some $a_1$, $a_2 > 0$). Hence we obtain
\begin{equation}\label{finiteLevyMeasureInequality}
 \mathbb{E}|X_t(x) - Y_t(y)| \leq \widetilde{L} e^{-\theta t}(|x-y| + 1) \,,
\end{equation}
 for some $\widetilde{L} > 0$, which is, however, undesirable since in order to be able to apply Theorem \ref{transportationInequalitiesTheorem} we would like to have $|x-y|$ and not $|x-y| + 1$ on the right hand side (cf. Remark \ref{finiteLevyMeasureRemark}). Thus we need to improve on the result from \cite{majka} and get an inequality like (\ref{majkaThm12}) but with a continuous function $f$ (i.e., with $a = 0$ in (\ref{formulaF})). To this end, we define 
\begin{equation}\label{Cepsilon}
 C_{\varepsilon} := 2 \int_0^{ \varepsilon/4  } |y|^2 \nu^L_1(dy) \,,
\end{equation}
where $\nu^L_1$ is the first marginal of the rotationally invariant measure $\nu^L$. The choice of $\varepsilon / 4$ as the upper integration limit is motivated by the calculations in the proof of Theorem 1.1 in \cite{majka}, see also the proof of Theorem \ref{jumpCouplingMainTheorem} below. Now consider a new condition.

\begin{myassumption}{L5}\label{Assumption6}
(Sufficient concentration of $\nu^L$ around zero)
For any $\lambda > 0$ there exists a $K(\lambda) > 0$ such that for all $\varepsilon < \lambda$ we have $\varepsilon \leq  K(\lambda) C_{\varepsilon}$. In other words, $\varepsilon / C_{\varepsilon}$ is bounded near zero or, using the big $O$ notation, $\varepsilon = O(C_{\varepsilon})$ as $\varepsilon \to 0$.
\end{myassumption}

Intuitively, it is an assumption about sufficient concentration of the L\'{e}vy measure $\nu^L$~around zero (sufficient small jump activity). It is satisfied e.g. for $\alpha$-stable processes with $\alpha \in [1,2)$ since in this case $C_\varepsilon = A \varepsilon^{2 - \alpha}$ for some constant $A = A(\alpha)$ and we have $\varepsilon / C_{\varepsilon} = A\varepsilon^{\alpha - 1}$.

It turns out that once we replace Assumptions \ref{Assumption3} and \ref{Assumption4} in Theorem 1.1 in \cite{majka} with Assumption \ref{Assumption6}, we are able to obtain (\ref{majkaThm12}) with a continuous function $f$, which is exactly what we need for Theorem \ref{transportationInequalitiesTheorem}. This is done in Section \ref{sectionCouplingPureJump} in Theorem \ref{jumpCouplingMainTheorem}. However, we are able to generalize this result even further.

\begin{theorem}\label{mainCouplingTheorem}
 Consider an SDE of the form 
 \begin{equation}\label{twoNoisesSDE}
 dX_t = b(X_t)dt + \sigma_1 dB^1_t + \sigma(X_t) dB^2_t + dL_t + \int_{U} g(X_{t-},u) \widetilde{N}(dt,du) \,,
\end{equation}
where $(B_t^1)_{t \geq 0}$ and $(B_t^2)_{t \geq 0}$ are $d$-dimensional Brownian motions, $(L_t)_{t \geq 0}$ is a pure jump L\'{e}vy process with L\'{e}vy measure $\nu^L$ satisfying Assumptions \ref{Assumption1}-\ref{Assumption2} and \ref{Assumption6}, whereas $\widetilde{N}$ is a compensated Poisson random measure on $\mathbb{R}_{+} \times U$ with intensity measure $dt \, \nu(du)$. Assume that all the sources of noise are independent, $\sigma_1 \in \mathbb{R}^{d \times d}$ is a constant matrix and the coefficients $b: \mathbb{R}^d \to \mathbb{R}^d$, $\sigma: \mathbb{R}^d \to \mathbb{R}^{d \times d}$ and $g: \mathbb{R}^d \times U \to \mathbb{R}^d$ satisfy Assumption \ref{AssumptionD1}. If at least one of the following two conditions is satisfied
\begin{enumerate}
 \item $\det \sigma_1 > 0$,
 \item $L_t \neq 0$ and Assumption \ref{AssumptionD2},
\end{enumerate}
then there exists a coupling $(X_t,Y_t)_{t \geq 0}$ of solutions to (\ref{twoNoisesSDE}) and constants $\widetilde{C}$, $\widetilde{c} > 0$ such that for any $x$, $y \in \mathbb{R}^d$ and any $t > 0$ we have
\begin{equation}\label{initialPerturbationInequality}
\mathbb{E}|X_t(x) - Y_t(y)| \leq \widetilde{C} e^{-\widetilde{c}t} |x-y| \,.
\end{equation}
\end{theorem}

 \begin{remark}
  The reason for the particular form of the equation (\ref{twoNoisesSDE}) is that in order to construct a coupling leading to the inequality (\ref{initialPerturbationInequality}) we need a suitable additive component of the noise. We can either use $(B_t^1)_{t \geq 0}$ if the condition (1) holds, or $(L_t)_{t \geq 0}$ if the condition (2) holds. The constants $\widetilde{C}$ and $\widetilde{c}$ depend on which noise we use. In particular, the constant $\widetilde{c}$ is either equal to $c$ defined by (\ref{defcEberle}) if we use $(B_t^1)_{t \geq 0}$ or to $c_1$ defined by (\ref{defc1}) if we use $(L_t)_{t \geq 0}$. On the other hand, if we have only a multiplicative Gaussian noise but the coefficient $\sigma$~is such that $\sigma \sigma^T$ is uniformly positive definite, we can use Lemma \ref{diffusionSplittingLemma} below to decompose this noise and extract an additive component satisfying (1). Without such an assumption on $\sigma$, Remark 2 in \cite{eberle} indicates that it might still be possible to perform a suitable construction, using the so-called Kendall-Cranston coupling, although this might significantly increase the level of sophistication of the proof. In the case of the jump noise, as far as we know there are currently no methods for obtaining couplings leading to inequalities like (\ref{initialPerturbationInequality}) in the case of purely multiplicative noise, and the recent papers treating this kind of problems (see e.g. \cite{jwang4}, \cite{majka} and \cite{luowang}) use methods that rely on the noise having at least some additive component. 
 \end{remark}

\begin{remark}\label{weakSolutionRemark}
 The coupling process $(X_t,Y_t)_{t \geq 0}$ is constructed as a unique strong solution to some $2d$-dimensional SDE. This allows us to infer that $(X_t,Y_t)_{t \geq 0}$ is in fact a~Markov process (see e.g. Theorem 6.4.5 in \cite{apple} or Proposition 4.2 in \cite{albeverio}, where it is shown how the Markov property follows from the uniqueness in law of solutions to SDEs with jumps). As a~consequence, we see that the inequality (\ref{initialPerturbationInequality}) actually implies that for any $0 \leq s \leq t$ we have
 \begin{equation*}
  \mathbb{E}[|X_t - Y_t| / \mathcal{F}_s ] \leq \widetilde{C} e^{-\widetilde{c}(t-s)} |X_s-Y_s| \,.
 \end{equation*}
\end{remark}
 
 \begin{remark}\label{finiteLevyMeasureRemark}
  Theorem \ref{mainCouplingTheorem} is obtained based on Theorem \ref{jumpCouplingMainTheorem} which is presented later in this paper. It is however possible to obtain analogous (but perhaps less useful) result based on the already mentioned Theorem 1.1 in \cite{majka}, where we have Assumptions \ref{Assumption3} and \ref{Assumption4} instead of Assumption \ref{Assumption6}. Then we get an inequality of the form (\ref{finiteLevyMeasureInequality}). It is still possible to obtain some transportation inequalities if in Theorem \ref{transportationInequalitiesTheorem} we replace the condition (\ref{initialChangeBound}) with a condition like (\ref{finiteLevyMeasureInequality}), but because of its form it forces us to additionally assume that the underlying intensity measure is finite (see Remark \ref{finiteIntensityRemark}). 
 \end{remark}

The above result is proved using the coupling methods developed in \cite{majka} and \cite{eberle}, and is of independent interest, as it extends some of the results obtained there. In particular, it immediately allows us to obtain exponential (weak) contractivity of the transition semigroup $(p_t)_{t \geq 0}$ associated with the solution to (\ref{twoNoisesSDE}), with respect to the $L^1$-Wasserstein distance $W_1$, as shown by the following corollary.

\begin{corollary}\label{corollaryExponentialContractivity}
 Under the assumptions of Theorem \ref{mainCouplingTheorem},
 \begin{equation*}
  W_1(\eta p_t, \mu p_t) \leq \widetilde{C}e^{-\widetilde{c}t} W_1(\eta, \mu)
 \end{equation*}
for any probability measures $\eta$ and $\mu$ on $\mathbb{R}^d$ and for any $t >0$. Moreover, $(p_t)_{t \geq 0}$ has an invariant measure $\mu_0$ and we have
\begin{equation*}
 W_1(\eta p_t, \mu_0) \leq \widetilde{C}e^{-\widetilde{c}t} W_1(\eta, \mu_0)
\end{equation*}
for any probability measure $\eta$ on $\mathbb{R}^d$ and any $t > 0$.
\end{corollary}

 This result follows immediately from (\ref{initialPerturbationInequality}) like in the proof of Corollary 3 in \cite{eberle} or the beginning of Section 3 in \cite{komorowski}. Using couplings allows us also to prove a related result involving a perturbation of the solution to (\ref{twoNoisesSDE}) by a change in the drift. This gives us a tool to determine some concrete cases in which the assumption (\ref{driftChangeBound}) from Theorem \ref{transportationInequalitiesTheorem} holds.

\begin{theorem}\label{theoremDriftPerturbation}
  Let $(X_t)_{t \geq 0}$ be like in Theorem \ref{mainCouplingTheorem} and suppose additionally that Assumption \ref{AssumptionD2} holds, $\det \sigma_1 > 0$ and that the coefficients $\sigma$ and $g$ are Lipschitz. Consider a process $(\widetilde{X}_t)_{t \geq 0}$ which is a~solution to (\ref{twoNoisesSDE}) with the drift perturbed by $u_t$, i.e.,
\begin{equation*}
 d\widetilde{X}_t = b(\widetilde{X}_t)dt + u_t dt + \sigma_1 dB^1_t + \sigma(\widetilde{X}_t) dB^2_t + dL_t + \int_{U} g(\widetilde{X}_{t-},u) \widetilde{N}(dt,du) \,,
\end{equation*}
where $u_t$ is either $\sigma_1 h_t$ or $\sigma(\widetilde{X}_t) h_t$ for some adapted $d$-dimensional process $h_t$. Then there exists a process $(Y_t)_{t \geq 0}$ such that $(X_t, Y_t)_{t \geq 0}$ is a coupling of solutions to (\ref{twoNoisesSDE}) and for any $0 \leq s \leq t$ we have
\begin{equation}\label{driftPerturbationInequality}
 \mathbb{E} [|\widetilde{X}_t - Y_t| / \mathcal{F}_s ] \leq C \int_s^t e^{c(r-(t-s))} |u_r| dr \,,
\end{equation}
where the constants $C$, $c > 0$ are given by (\ref{defCEberle}) and (\ref{defcEberle}), respectively.
\end{theorem}

Observe that the constants above depend on the function $\kappa$ and hence to calculate their explicit values we need to apply the right version of $\kappa$ in the formulas (\ref{defCEberle}) and (\ref{defcEberle}), i.e., the version that is used in the proof of Theorem \ref{theoremDriftPerturbation}. Now, combining Theorems \ref{mainCouplingTheorem} and \ref{theoremDriftPerturbation} to check validity of assumptions of Theorem \ref{transportationInequalitiesTheorem}, we get the following result.
 
 \begin{corollary}\label{corollaryTransportation}
  Consider the setup of Theorem \ref{mainCouplingTheorem}. Suppose all its assumptions and Assumption \ref{AssumptionD2} are satisfied and additionally that $\det \sigma_1 > 0$ and the coefficients $\sigma$ and $g$ are Lipschitz. Moreover, assume that $(X_t)_{t \geq 0}$ is Malliavin differentiable ($X_t \in \mathbb{D}^{1,2}$ for all $t \geq 0$) and, similarly to Theorem \ref{transportationInequalitiesTheorem}, that there exists a constant $\sigma_{\infty}$ such that for any $x \in \mathbb{R}^d$ we have $\| \sigma(x) \| \leq \sigma_{\infty}$ and there exists a measurable function $g_{\infty}: U \to \mathbb{R}$ such that $|g(x,u)| \leq g_{\infty}(u)$ for any $x \in \mathbb{R}^d$ and $u \in U$. Assume further that there exists some $\lambda > 0$ such that Assumption \ref{AssumptionE} is satisfied and that there exists $\widetilde{\lambda} > 0$ such that
  \begin{equation*}
       \beta^L(\widetilde{\lambda}) := \int_{U} (e^{\widetilde{\lambda} u } - \widetilde{\lambda} u - 1) \nu^L(du) < \infty \,.
  \end{equation*}
  Then the transportation inequality (\ref{transIneq1}) from the statement of Theorem \ref{transportationInequalitiesTheorem} holds with 
  \begin{equation*}
  \begin{split}
   \alpha_{T}(r) &:= \sup_{\lambda \geq 0} \Bigg\{ r \lambda - \int_0^T \beta(\widetilde{C}e^{-\widetilde{c}(T-t)} \lambda)dt - \int_0^T \beta^L(\widetilde{C}e^{-\widetilde{c}(T-t)} \lambda)dt \\
   &-\frac{(\sigma_{\infty}^2 + \| \sigma_1 \|^2) \lambda^2}{2} C^2 \frac{1-e^{-2cT}}{2c} \Bigg\} \,.
   \end{split}
  \end{equation*}
Moreover, for the invariant measure $\mu_0$ we have
\begin{equation}\label{alphaW1HforInvariant}
 \alpha_{\infty}(W_1(\eta, \mu_0)) \leq H(\eta | \mu_0)
\end{equation}
for any probability measure $\eta$ on $\mathbb{R}^d$, with $\alpha_{\infty}$ defined as the pointwise limit of $\alpha_T$ as $T \to \infty$. Finally, the inequality (\ref{transIneq2}) holds with
\begin{equation*}
\begin{split}
 \alpha^P_T(r) &:= \sup_{\lambda \geq 0} \Bigg\{ r \lambda - \int_0^T \beta\left(\lambda \widetilde{C} \frac{1-e^{-\widetilde{c}(T-t)}}{\widetilde{c}}  \right)dt - \int_0^T \beta^L\left(\lambda \widetilde{C} \frac{1-e^{-\widetilde{c}(T-t)}}{\widetilde{c}}  \right)dt\\
 &- \frac{(\sigma_{\infty}^2 + \| \sigma_1 \|^2) \lambda^2}{2} C^2 \int_0^T \left( \frac{1-e^{-c(T-t)}}{c} \right)^2 dt \Bigg\} \,.
 \end{split}
\end{equation*} 
The constants $\widetilde{c}$, $\widetilde{C}$, $c$ and $C$ appearing in the definitions of $\alpha_T$ and $\alpha_T^P$ are the same as in (\ref{initialPerturbationInequality}) and (\ref{driftPerturbationInequality}).
 \end{corollary}

 This corollary extends the results from Theorem 2.2 in \cite{lwu} to the case where we drop the global dissipativity assumption required therein, as long as we have an additive component of the noise, which we can use in order to construct a coupling required in our method. It is easy to notice that the corollaries in Section 2 in \cite{lwu} (various results regarding concentration of measure for solutions of (\ref{twoNoisesSDE}) in the pure jump case) hold as well under our assumptions. We also extend Theorem 2.2 from \cite{yma}, where similar results are proved in the jump diffusion case under assumptions analogous to the ones in \cite{lwu}. However, in \cite{yma} there are additionally stronger assumptions on regularity of the coefficients, which are needed to get bounds on Malliavin derivatives of solutions to (\ref{twoNoisesSDE}). Here we use a different method of getting such bounds (cf. Remark \ref{remarkMalliavinBounds}) which does not require coefficients to be differentiable and works whenever we have $X_t \in \mathbb{D}^{1,2}$ for all $t \geq 0$.

 \begin{example}
  To have a jump noise satisfying all the assumptions of Corollary \ref{corollaryTransportation}, we can take a L\'{e}vy process whose L\'{e}vy measure behaves near the origin like that of an $\alpha$-stable process with $\alpha \in (1,2)$ (so that Assumptions \ref{Assumption1}-\ref{Assumption2} and \ref{Assumption6} are satisfied), but has exponential moments as well (so that Assumption \ref{AssumptionE} is also satisfied). A natural example of such a process is the so-called relativistic $\alpha$-stable process, which is a L\'{e}vy process $(L_t)_{t \geq 0}$ with the characteristic function given by
  \begin{equation*}
   \mathbb{E}\exp \left( i \langle z , L_t \rangle \right) = \exp \left(-t\left[ (m^{1/\beta} + |z|^2)^{\beta} - m \right] \right)
  \end{equation*}
for $z \in \mathbb{R}^d$, with $\beta = \alpha / 2$ and some parameter $m > 0$. For more information on this process, see e.g. \cite{carmona} where Corollary II.2 and Proposition II.5 show that it indeed satisfies Assumption \ref{AssumptionE}, or \cite{ryznar} where in Lemma 2 the formula for the density of its L\'{e}vy measure is calculated, from which we can easily see that Assumption \ref{Assumption6} holds. SDEs driven by relativistic stable processes (and in fact also by a significantly more general type of noise) have been recently studied in \cite{huang}.
 \end{example}

 \begin{remark}
  Both in \cite{lwu} and \cite{yma}, apart from the transportation inequalities for measures $\delta_x p_T$ on $\mathbb{R}^d$ and for measures $\mathbb{P}_{x,[0,T]}$ on the path space with the $L^1$ metric, there are also inequalities on the path space with the $L^{\infty}$ metric defined by $d_{\infty}(\gamma_1, \gamma_2) = \sup_{t \in [0,T]} |\gamma_1(t) - \gamma_2(t)|$ (see Theorem 2.11 in \cite{lwu} and Theorem 2.8 in \cite{yma}). However, the method of proof for these (see the second part of the proof of Lemma 3.3 in \cite{yma}) involves proving an inequality of the type
  \begin{equation*}
   \mathbb{E} \sup_{0 \leq s \leq t} |X_s(x) - X_s(y)|^2 \leq e^{\widehat{C}t}|x-y|^2 \,,
  \end{equation*}
with some constant $\widehat{C}$, which requires the integral form of the Gronwall inequality, which can only work if the constant $\widehat{C}$ is positive (cf. Remark 2.3 in \cite{shao}). Since this is the case even under the global dissipativity assumption, we have not been able to use couplings to improve on these results in any way and hence we skip them in our presentation, referring the interested reader to \cite{lwu} and \cite{yma}.
 \end{remark}

 \begin{remark}\label{remarkShao}
  Another possible application of our approach would be to extend the results from \cite{shao}, where transportation inequalities were studied in the context of regime switching processes, modelled by stochastic differential equations with both Gaussian and Poissonian noise (see (2.1) and (2.2) therein). There a kind of one-sided Lipschitz condition is imposed on the coefficients (see the condition (A3) in \cite{shao}) and, as pointed out in Remark 2.2 therein, transportation inequalities on the path space can be obtained without dissipativity. However, in such a case the constants with which those inequalities hold for $\mathbb{P}_{x,[0,T]}$, explode when $T \to \infty$ (see Theorem 2.1 in \cite{shao}). Since the method of proof used in \cite{shao} is a direct extension of the one developed by Liming Wu in \cite{lwu}, it should be possible to apply our reasoning to obtain non-exploding constants at least in (2.13) in \cite{shao} under a dissipativity at infinity condition. This is, however, beyond the scope of our present paper.
 \end{remark}

 \begin{remark}\label{remarkOtherPapers}
  In the present paper, we only explain in detail how to check assumptions of Theorem \ref{transportationInequalitiesTheorem} using the approach of Theorems \ref{mainCouplingTheorem} and \ref{theoremDriftPerturbation}. However, it may be possible to obtain inequalities like (\ref{initialChangeBound}) and (\ref{driftChangeBound}) by other methods. For example, in a~recent paper \cite{luowang}, D. Luo and J. Wang obtained an inequality like (\ref{initialChangeBound}) for equations of the type $dX_t = b(X_t)dt + dL_t$ under different than ours assumptions on the L\'{e}vy measure and using a different coupling (see Theorem 1.1 therein; (1.6) in \cite{luowang} follows from an inequality like (\ref{initialChangeBound}) which is needed in the proof of Theorem 3.1 therein). This is sufficient to get the transportation inequalities like in our Theorem \ref{transportationInequalitiesTheorem} for an SDE with pure jump noise under their set of assumptions (plus, additionally, Assumption \ref{AssumptionE}). On the other hand, Eberle, Guillin and Zimmer in \cite{eberleguillin} showed an inequality like (\ref{initialChangeBound}) for equations of the type $dX_t = b(X_t)dt + dB_t$ without assuming dissipativity even at infinity, at the cost of multiplying the right hand side of (\ref{initialChangeBound}) by a factor which, however, can possibly be controlled under some suitable integrability assumptions for $(X_t)_{t \geq 0}$, cf. Theorem 2 and formula (28) in \cite{eberleguillin}. This could lead to obtaining at least some concentration inequalities like (\ref{concentrationIneq1}) for solutions of equations of the type $dX_t = b(X_t)dt + dB_t + \int_U g(X_{t-},u) \widetilde{N}(dt,du)$, under some weaker than ours assumptions on the coefficients $b$ and $g$. These examples show robustness of our formulation of Theorem \ref{transportationInequalitiesTheorem}, as it allows us to easily obtain transportation or concentration inequalities in many cases where inequalities like (\ref{initialChangeBound}) arise naturally.  
 \end{remark}

 The crucial step in the proof of Theorem \ref{transportationInequalitiesTheorem} is to find upper bounds for the Malliavin derivatives of $X_t$. Thus, in the process of proving our main results, we also obtain some bounds that might be interesting on their own in the context of the Malliavin calculus.

 \begin{theorem}\label{generalMalliavinResult}
  Let $(X_t)_{t \geq 0}$ be a Malliavin differentiable solution to (\ref{generalSDEsect2}) such that there exists a coupling $(X_t, Y_t')_{t \geq 0}$ for which (\ref{driftChangeBound}) holds. Assume that there exists a constant $\sigma_{\infty}$ such that for any $x \in \mathbb{R}^d$ we have $\| \sigma(x) \| \leq \sigma_{\infty}$. Then for any Lipschitz functional $f : \mathbb{R}^d \to \mathbb{R}$ with $\| f \|_{\operatorname{Lip}} \leq 1$, for any adapted, $\mathbb{R}_{+}$-valued process $g$ and for any $0 \leq s \leq r \leq t$ we have
  \begin{equation}\label{MalliavinExpectationEstimate}
  \mathbb{E} \int_s^r g_u |\mathbb{E}[ \nabla_u f(X_t) | \mathcal{F}_u]|^2 du \leq c^2_2(t) \sigma_{\infty} \mathbb{E} \int_s^r g_u c^2_3(u) du \,.
  \end{equation}
  Moreover, we have
  \begin{equation}\label{MalliavinLinftyEstimate}
  \| \mathbb{E}[ \nabla_u f(X_t) | \mathcal{F}_u] \|_{L^{\infty}(\Omega \times [0,t])} \leq c_2(t) \sigma_{\infty} \sup_{u \leq t} c_3(u) \,,
  \end{equation}
where the $L^{\infty}$ norm is the essential supremum on $\Omega \times [0,t]$.
 \end{theorem}

 On the other hand, using the condition (\ref{initialChangeBound}), we can obtain related bounds for the Malliavin derivative $D$ of Lipschitz functionals of $X_t$ with respect to the Poisson random measure $N$ (see Section \ref{MalliavinPoissonian} for details).
 
 In the same way in which Corollary \ref{corollaryTransportation} follows from Theorem \ref{transportationInequalitiesTheorem} via Theorems \ref{mainCouplingTheorem} and \ref{theoremDriftPerturbation}, the following corollary follows from Theorem \ref{generalMalliavinResult} via Theorem \ref{theoremDriftPerturbation}.
 
 \begin{corollary}\label{corollaryMalliavinBounds}
  Let $(X_t)_{t \geq 0}$ be a Malliavin differentiable solution to (\ref{twoNoisesSDE}), satisfying the assumptions of Theorem \ref{mainCouplingTheorem} with $\det \sigma_1 > 0$ and $\lim_{r \to 0} r \kappa(r) = 0$ (i.e., Assumption \ref{AssumptionD2}). Moreover, assume that the coefficients $\sigma$ and $g$ are Lipschitz and that there exists a~constant $\sigma_{\infty}$ such that for any $x \in \mathbb{R}^d$ we have $\| \sigma(x) \| \leq \sigma_{\infty}$. Denote by $\nabla^i$ the Malliavin derivative with respect to $(B_t^i)_{t \geq 0}$ for $i \in \{ 1, 2\}$. Then for any functional $f$ and any process $g$ like above and for any $0 \leq s \leq r \leq t$ we have
\begin{equation}\label{additiveMalliavinEstimate}
 \mathbb{E} \int_s^r g_u |\mathbb{E}[ \nabla_u^1 f(X_t) | \mathcal{F}_u]|^2 du \leq C^2 \| \sigma_1 \|^2 \mathbb{E} \int_s^r g_u  e^{2c(u-t)} du
\end{equation}
and
\begin{equation}\label{multiplicativeMalliavinEstimate}
 \mathbb{E} \int_s^r g_u |\mathbb{E}[ \nabla_u^2 f(X_t) | \mathcal{F}_u]|^2 du \leq C^2  \sigma_{\infty}^2  \mathbb{E} \int_s^r g_u e^{2c(u-t)} du \,,
\end{equation}
where $C$ and $c$ are the same as in (\ref{driftPerturbationInequality}). We also have $L^{\infty}$ bounds analogous to (\ref{MalliavinLinftyEstimate}) for $\nabla_u^1 f(X_t)$ and $\nabla_u^2 f(X_t)$, with the upper bound being, respectively, $\|\sigma_1 \|$ and $\sigma_{\infty}$.
 \end{corollary}

In analogy to our comment below the statement of Theorem \ref{generalMalliavinResult}, we observe here that a related corollary for the Malliavin derivatives with respect to $(L_t)_{t \geq 0}$ and $N$ is also true (see the end of Section \ref{MalliavinPoissonian}, in particular (\ref{twoNoisesMalliavinPoissonBound1}) and (\ref{twoNoisesMalliavinPoissonBound2})).
 
 \begin{remark}\label{remarkMalliavinBounds}
  If the global dissipativity assumption is satisfied and the coefficients in the equation are continuously differentiable, it is possible to obtain much stronger bounds than (\ref{additiveMalliavinEstimate}) and (\ref{multiplicativeMalliavinEstimate}). Namely, for the multiplicative noise we get
  \begin{equation*}
   \mathbb{E}[ \| \nabla_s X_t \|^2_{HS} | \mathcal{F}_s] \leq \| \sigma(X_s) \|^2_{HS} e^{2K(s-t)}
  \end{equation*}
for any $0 \leq s \leq t$, where $K > 0$ is the constant with which the global dissipativity condition holds (see Lemma 3.4 in \cite{yma}). We were not able to obtain such bounds in our case. However, our assumptions are much weaker than the ones in \cite{yma} and the bounds (\ref{additiveMalliavinEstimate}) and (\ref{multiplicativeMalliavinEstimate}) are sufficient to prove the transportation inequalities in Corollary \ref{corollaryTransportation}. On the other hand, our bounds for the Malliavin derivative $D$ with respect to the Poisson random measure $N$ have the same form as the ones in \cite{lwu} and \cite{yma} (cf. Section \ref{MalliavinPoissonian} in the present paper and Section 4.2 in \cite{lwu}).
 \end{remark}

 The remainder of the paper is organized as follows. In Section \ref{sectionCouplingPureJump} we present an extension of the results from \cite{majka} regarding couplings of solutions to SDEs driven by pure jump L\'{e}vy noise. In Section \ref{sectionCouplingJumpDiffusion} we explain how to further extend these results to the case of more general jump diffusions and hence we prove Theorem \ref{mainCouplingTheorem}. In Section \ref{subsectionBrownianMalliavin} we introduce our technique of obtaining estimates like (\ref{driftChangeBound}) in Lemma \ref{lemmaForMalliavin}, which then leads directly to the proofs of Theorem \ref{theoremDriftPerturbation} and Theorem \ref{generalMalliavinResult}, followed by the proof of Corollary \ref{corollaryMalliavinBounds}. In Section \ref{MalliavinPoissonian} we explain how to show related results in the case of Malliavin derivatives with respect to Poisson random measures. In Section \ref{sectionProofsTransportation} we finally prove the transportation and concentration inequalities, i.e., Theorem \ref{transportationInequalitiesTheorem}, Theorem \ref{concentrationInequalitiesTheorem} and Corollary \ref{corollaryTransportation}.

\section{Coupling of SDEs with pure jump noise}\label{sectionCouplingPureJump}
 
Here we consider an SDE of the form
\begin{equation}\label{purejumpSDE}
 dX_t = b(X_t)dt + dL_t \,,
\end{equation}
where $(L_t)_{t \geq 0}$ is a pure jump L\'{e}vy process and the drift function $b$ is continuous and satisfies a one-sided Lipschitz condition. In this section, let $N$ be the Poisson random measure on $\mathbb{R}_{+} \times \mathbb{R}^d$ associated with $(L_t)_{t \geq 0}$ via
\begin{equation*}
 L_t = \int_0^t \int_{\{|v|>1\}} v N(ds,dv) + \int_0^t \int_{\{|v| \leq 1\}} v \widetilde{N}(ds,dv)
\end{equation*}
and let $dt \, \nu(dv)$ be its intensity measure. Following Section 2.2 in \cite{majka}, we can replace $N$ with a Poisson random measure on $\mathbb{R}_{+} \times \mathbb{R}^d \times [0,1]$ with intensity $dt \, \nu(dv) \, du$, where $du$ is the Lebesgue measure on $[0,1]$, thus introducing an additional control variable $u \in [0,1]$. By a slight abuse of notation, we keep denoting this new Poisson random measure by $N$.
 We can thus write (\ref{purejumpSDE}) as
 \begin{equation*}
 dX_t = b(X_t)dt + \int_{\{ |v| > 1 \} \times [0,1]} v N(dt,dv,du) + \int_{\{ |v| \leq 1 \} \times [0,1]} v \widetilde{N}(dt,dv,du) \,.
\end{equation*}
Without loss of generality, we can choose a constant $m > 1$ and rewrite the equation above as
 \begin{equation}\label{eqX}
 dX_t = b(X_t)dt + \int_{\{ |v| > m \} \times [0,1]} v N(dt,dv,du) + \int_{\{ |v| \leq m \} \times [0,1]} v \widetilde{N}(dt,dv,du) \,.
\end{equation}
Formally we should then change the drift function by an appropriate constant, but since such an operation does not change any relevant properties of the drift, we choose to keep denoting the drift by $b$.
Now we can define a coupling $(X_t,Y_t)_{t \geq 0}$ by putting
\begin{equation}\label{optimalY}
 \begin{split}
  dY_t &= b(Y_t)dt + \int_{\{ |v| > m \} \times [0,1]} v N(dt,dv,du) \\
  &+ \int_{\{ |v| \leq m \} \times [0,1]} (X_{t-} - Y_{t-} + v)\mathbf{1}_{\{ u < \rho(v, Z_{t-}) \}} \widetilde{N}(dt,dv,du) \\
  &+ \int_{\{ |v| \leq m \} \times [0,1]} R(X_{t-},Y_{t-})v\mathbf{1}_{\{ u \geq \rho(v, Z_{t-}) \}} \widetilde{N}(dt,dv,du) \,,
 \end{split}
\end{equation}
for $t < T := \inf \{ t > 0 : X_t = Y_t \}$ and $Y_t = X_t$ for $t \geq T$, where $Z_t := X_t - Y_t$,
\begin{equation*}
  \rho(v, Z_{t-}) := \frac{q(v) \wedge q(v + Z_{t-})\mathbf{1}_{\{ |v+Z_{t-}| \leq m \}}}{q(v)}
\end{equation*}
if $q(v) \neq 0$ and $\rho(v, Z_{t-}) := 1$ if $q(v) = 0$, where $\nu(dv) = q(v)dv$ and
\begin{equation}\label{reflectionOperatorLevy}
 R(X_{t-},Y_{t-}) := I - 2\frac{(X_{t-}-Y_{t-})(X_{t-}-Y_{t-})^T}{|X_{t-} - Y_{t-}|^2} = I - 2e_{t-}e_{t-}^T \,,
\end{equation}
with $e_t := (X_t - Y_t)/|X_t- Y_t|$. Intuitively, it is a combination of a modification of the reflection coupling with a positive probability of bringing the marginal processes together instead of performing the reflection (for jumps of size smaller than $m$) and the synchronous coupling (for jumps larger than $m$). We can call it the mirror coupling. For the coupling construction itself, $m$ can be chosen arbitrarily. For obtaining convergence rates in Wasserstein distances, we choose $m$ based on Assumption \ref{Assumption3} satisfied by the L\'{e}vy measure $\nu$ of $(L_t)_{t \geq 0}$. For the discussion explaining this construction in detail see Section 2 in \cite{majka}.

Under Assumptions \ref{Assumption1} and \ref{Assumption2} it has been proved in \cite{majka} (see Theorem 1.1 therein) that the $2d$-dimensional SDE given by (\ref{eqX}) and (\ref{optimalY}) has a unique strong solution which is a coupling of solutions to (\ref{purejumpSDE}). Then this coupling was used to prove that, under additional Assumptions \ref{Assumption3} and \ref{Assumption4} and a dissipativity at infinity condition on the drift, the inequality (\ref{majkaThm12}) holds with a discontinuous function $f$, i.e., we have
\begin{equation*}
 \mathbb{E}|X_t(x) - Y_t(y)| \leq L e^{-\theta t}f(|x-y|)
\end{equation*}
for some constants $L > 1$ and $\theta > 0$.

Now we turn to the proof of a modification of the main result in \cite{majka}, which will give us an inequality like (\ref{majkaThm12}), but with a continuous function $f$. Recall that in the case of an equation of the form (\ref{purejumpSDE}), the function $\kappa$ is such that for all $x$, $y \in \mathbb{R}^d$ we have
\begin{equation}\label{defKappa}
 \langle b(x) - b(y) , x - y \rangle \leq - \kappa(|x-y|)|x-y|^2 \,.
\end{equation}

\begin{theorem}\label{jumpCouplingMainTheorem}
  Let $(X_t)_{t \geq 0}$ be a Markov process in $\mathbb{R}^d$ given as a solution to the stochastic differential equation (\ref{purejumpSDE}), where $(L_t)_{t \geq 0}$ is a pure jump L\'{e}vy process satisfying Assumptions \ref{Assumption1}-\ref{Assumption2} and Assumption \ref{Assumption6} and $b: \mathbb{R}^d \to \mathbb{R}^d$ is a continuous, one-sided Lipschitz vector field satisfying Assumptions \ref{AssumptionD1} and \ref{AssumptionD2}. Then there exists a coupling of solutions to (\ref{purejumpSDE}) defined as a strong solution to the $2d$-dimensional SDE given by (\ref{eqX}) and (\ref{optimalY}) and a continuous concave function $f_1 : \mathbb{R}_{+} \to \mathbb{R}_{+}$ such that
 \begin{equation*}
  \mathbb{E}f_1(|X_t(x) - Y_t(y)|) \leq e^{-c_1t}f_1(|x-y|)
 \end{equation*}
holds with some constant $c_1 > 0$ for any $t > 0$ and any $x$, $y \in \mathbb{R}^d$. By the construction of $f_1$, we also have
 \begin{equation*}
  \mathbb{E}|X_t(x) - Y_t(y)| \leq Le^{-c_1t}|x-y|
 \end{equation*}
with some constant $L > 0$.

\begin{proof}
 The existence of the coupling as a strong solution to the system (\ref{eqX})-(\ref{optimalY}) has been proved in Section 2 in \cite{majka}. Now we will explain how to modify the proof of the inequality (1.8) in Theorem 1.1 in \cite{majka} in order to prove the new result presented here. Denote
\begin{equation*}
 Z_t := X_t - Y_t \,.
\end{equation*}
Using the expression (\ref{optimalY}) for $dY_t$, we can write
\begin{equation*}
  \begin{split}
  dZ_t &= (b(X_t) - b(Y_t))dt + \int_{\{ |v| \leq m \} \times [0,1]} (I-R(X_{t-},Y_{t-}))v \widetilde{N}(dt,dv,du) \\
  &+ \int_{\{ |v| \leq m \} \times [0,1]} A(X_{t-},Y_{t-},v,u) \widetilde{N}(dt,dv,du) \,.
 \end{split}
\end{equation*}
where $A(X_{t-},Y_{t-},v,u):= -(Z_{t-} + v - R(X_{t-},Y_{t-})v)\mathbf{1}_{\{ u < \rho(v, Z_{t-}) \}}$.
Applying the It\^{o} formula (see e.g. Theorem 4.4.10 in \cite{apple}) with a function $f_1$ we get
\begin{equation}\label{ItoOptimal}
 \begin{split}
 f_1(|Z_t|) - f_1(|&Z_0|) = \int_0^t f_1'(|Z_{s-}|)\frac{1}{|Z_{s-}|}\langle Z_{s-} , b(X_{s-}) - b(Y_{s-}) \rangle ds \\
 &+\int_0^t \int_{\{ |v| \leq m \} \times [0,1]} f_1'(|Z_{s-}|)\frac{1}{|Z_{s-}|}\langle Z_{s-},(I-R(X_{s-},Y_{s-}))v \rangle \widetilde{N}(ds,dv,du)  \\
 &+\int_0^t \int_{\{ |v| \leq m \} \times [0,1]} f_1'(|Z_{s-}|)\frac{1}{|Z_{s-}|}\langle Z_{s-},A(X_{s-},Y_{s-},v,u) \rangle \widetilde{N}(ds,dv,du)\\
 &+\sum_{s \in (0,t]}  \left( |\Delta Z_s|^2 \int_0^1 (1-u) f_1''(|Z_{s-} + u \Delta Z_s|) du \right) \,.
 \end{split}
\end{equation}
where the last term is obtained from the usual sum over jumps appearing in the It\^{o} formula by applying the Taylor formula and using the fact that in our coupling the vectors $Z_{s-}$ and $\Delta Z_s$ are always parallel (see Section 3 in \cite{majka} for details). 
Now we introduce a sequence of stopping times $(\tau_n)_{n=1}^{\infty}$ defined by
\begin{equation}\label{stoppingTimes}
 \tau_n := \inf \{ t \geq 0 : |Z_t| \notin (1/n , n) \} \,.
\end{equation}
Note that we have $\tau_n \to T$ as $n \to \infty$, which follows from non-explosiveness of $(Z_t)_{t \geq 0}$. By some tedious but otherwise easy computations (see the proof of Theorem 1.1 in \cite{majka} for details, specifically Lemma 3.1 and Lemma 3.2 therein) we can show that 
\begin{equation}\label{majkaBound1}
 \mathbb{E} \int_0^{t \wedge \tau_n} \int_{\{ |v| \leq m \} \times [0,1]} f_1'(|Z_{s-}|)\frac{1}{|Z_{s-}|}\langle Z_{s-},(I-R(X_{s-},Y_{s-}))v \rangle \widetilde{N}(ds,dv,du) = 0 \,.
 \end{equation}
 and
  \begin{equation}\label{majkaBound2}
 \mathbb{E} \int_0^{t \wedge \tau_n} \int_{\{ |v| \leq m \} \times [0,1]} f_1'(|Z_{s-}|)\frac{1}{|Z_{s-}|} \langle Z_{s-},A(X_{s-},Y_{s-},v,u) \rangle \widetilde{N}(ds,dv,du) = 0 \,.
 \end{equation}
 In \cite{majka} it is also shown (see Lemma 3.3 therein) that for any $t > 0$, we have
 \begin{equation}\label{majkaBound3}
  \mathbb{E} \sum_{s \in (0,t]}  \left( |\Delta Z_s|^2 \int_0^1 (1-u) f_1''(|Z_{s-} + u \Delta Z_s|) du \right) \leq C_{\varepsilon} \mathbb{E} \int_0^t \bar{f}_{\varepsilon}(|Z_{s-}|) \mathbf{1}_{\{ |Z_{s-}| > \delta \}} ds \,,
 \end{equation}
where $0 < \delta < 2m$, $\varepsilon \leq \delta$, the constant $C_{\varepsilon}$ is defined as in (\ref{Cepsilon}) with the first marginal $\nu_1$ of the measure $\nu$ and the function $\bar{f}_{\varepsilon}$ is defined by
\begin{equation*}
 \bar{f}_{\varepsilon}(y):= \sup_{x \in (y - \varepsilon, y)}f_1''(x) \,.
\end{equation*}
It is important to note that in order for (\ref{majkaBound3}) to hold, $m$ has to be chosen in such a way that
\begin{equation*}
 \int_{ -\varepsilon/2 }^0 |y|^2 \nu^m_1(dy) \geq \int_{ -\varepsilon/4  }^0 |y|^2 \nu_1(dy) = \frac{C_{\varepsilon}}{2} \,,
\end{equation*}
where $\nu^m_1$ is the first marginal of the truncated measure $\nu^m(dv) := \mathbf{1}_{\{ |v| \leq m \}} \nu(dv)$. This is, however, not a problem, since $m$ can always be chosen large enough, cf. the discussion in Section 2.2 in \cite{majka}.
The crucial element of the proof in \cite{majka}, after getting the bounds (\ref{majkaBound1}), (\ref{majkaBound2}) and (\ref{majkaBound3}), is the construction of a function $f_1$ and a~constant $c_1 > 0$ such that
\begin{equation}\label{functionalInequality}
 -f_1'(r)\kappa(r)r + C_{\varepsilon} \bar{f}_{\varepsilon}(r) \leq -c_1 f_1(r)
\end{equation}
holds for all $r > \delta$, where $\kappa$ is the function satisfying (\ref{defKappa}). Combining this with (\ref{ItoOptimal}) and using Assumption \ref{Assumption3} and the discontinuity of the distance function to deal with the case of $r \leq \delta$ (see Lemma 3.7 in \cite{majka}), it is shown how to get a~bound of the form
\begin{equation*}
 \mathbb{E}f_1(|Z_{t \wedge \tau_n}|) - \mathbb{E}f_1(|Z_0|) \leq \mathbb{E} \int_0^{t \wedge \tau_n} -c_1 f_1(|Z_s|) ds \,,
\end{equation*}
which then leads to (\ref{majkaThm12}). Now we will show a different way of dealing with the case of $r \leq \delta$, using Assumption \ref{Assumption6} instead of Assumption \ref{Assumption3}, which allows us to keep the continuity of $f_1$.

It is quite easy to see (using once again the fact that $Z_{s-}$ and $\Delta Z_s$ are parallel, cf. the proof of Lemma 3.3 in \cite{majka}) that for any $u \in (0,1)$ we have
\begin{equation*}
 f_1''(|Z_{s-} + u \Delta Z_s|) \leq \sup_{x \in (|Z_{s-}|, |Z_{s-}| + \varepsilon)}f_1''(x)\mathbf{1}_{\{ |Z_s| \in (|Z_{s-}|,|Z_{s-}| + \varepsilon) \}} \,.
\end{equation*}
We also have 
\begin{equation*}
 \{ |Z_s| \in (|Z_{s-}|,|Z_{s-}| + \varepsilon) \} = \{ |Z_s| > |Z_{s-}|\} \cap \{ |\Delta Z_s| < \varepsilon \} \,,
\end{equation*}
and the condition $|Z_s| > |Z_{s-}|$ is equivalent to $\langle \Delta Z_s , 2Z_{s-}+ \Delta Z_s \rangle > 0$. Therefore, mimicking the argument in the proof of Lemma 3.3 in \cite{majka} we get that
 \begin{equation}\label{ItoJumpSum}
  \mathbb{E} \sum_{s \in (0,t]}  \left( |\Delta Z_s|^2 \int_0^1 (1-u) f_1''(|Z_{s-} + u \Delta Z_s|) du \right) \leq C_{\varepsilon} \mathbb{E} \int_0^t \hat{f}_{\varepsilon}(|Z_{s-}|) ds \,,
 \end{equation}
where
\begin{equation*}
 \hat{f}_{\varepsilon}(y):= \sup_{x \in (y - \varepsilon, y)}f_1''(x) \mathbf{1}_{\{ y > \delta \}} + \sup_{x \in (y, y + \varepsilon)}f_1''(x)\mathbf{1}_{\{ y \leq \delta \}} \,.
\end{equation*}

Now we will show that under Assumption \ref{Assumption6}, after a small modification in the formulas from \cite{majka}, the inequality
\begin{equation}\label{newFunctionalInequality}
 -f_1'(r)\kappa(r)r + C_{\varepsilon} \hat{f}_{\varepsilon}(r) \leq -c_1 f_1(r)
\end{equation}
holds for all $r > 0$ (note that here we have $\hat{f}_{\varepsilon}(r)$ in place of $\bar{f}_{\varepsilon}(r)$ in (\ref{functionalInequality})).
The function $f_1$, constructed in Lemma 3.6 in \cite{majka} in order to satisfy (\ref{functionalInequality}), is such that $f_1' \geq 0$, $f_1'' \leq 0$ and is defined in the following way 
\begin{equation}\label{defFMajka}
 f_1(r) = \int_0^r \phi(s)g(s) ds \,,
\end{equation}
where
\begin{align*}
 \phi(r) &:= \exp{\left( -\int_0^r \frac{\bar{h}(t)}{C_{\varepsilon}} dt \right)} \,, &\quad  g(r) &:= \begin{cases} 
      1 - \frac{c_1}{C_{\varepsilon}} \int_0^r \frac{\Phi(t+\varepsilon)}{\phi(t)} dt \,, & r\leq R_1 \,, \\
      \frac{1}{2} \,, & r\geq R_1 \,.
   \end{cases} 
\end{align*}
Here $R_1 > 0$ is given by formulas
\begin{equation}\label{Rzero}
\begin{split}
 R_0 &= \inf \left\{R \geq 0 : \forall r \geq R : \kappa(r) \geq \frac{2M}{R} \right\} \,, \\
 R_1 &= \inf \left\{ R \geq R_0 + \varepsilon : \forall r \geq R : \kappa(r) \geq \frac{2C_{\varepsilon}}{(R-R_0)R} + \frac{2M}{R} \right\} \,,
 \end{split}
 \end{equation}
 but can be chosen arbitrarily large if necessary and $c_1$ is a positive constant given by
\begin{equation}\label{defc1}
 c_1 := \frac{C_{\varepsilon}}{2} \left( \int_0^{R_1} \frac{\Phi(t + \varepsilon)}{\phi(t)} dt \right)^{-1}
\end{equation}
(cf. (3.29) in \cite{majka}). Moreover, we have
\begin{align*}
 \bar{h}(r)&:=\sup_{t \in (r,r+\varepsilon)}h^{-}(t) \,, &\enskip  \Phi(r) &:= \int_0^r \phi(s) ds 
\end{align*}
and $h^{-} = - \min \{ h, 0 \}$ is the negative part of the function
\begin{equation}\label{defH}
 h(r) := r\kappa(r) - 2M\,,
\end{equation}
with some $M > 0$ to be chosen later. Actually, in Lemma 3.6 in \cite{majka} the function $h$ is given by $h(r) := r\kappa(r)$, whereas $R_0$ and $R_1$ are chosen with $M = 0$ and this already gives (\ref{functionalInequality}). However, it is easy to check that by taking $M > 0$ we get 
\begin{equation}\label{functionalInequalityWithM}
 -f_1'(r)\kappa(r)r + 2f_1'(r)M + C_{\varepsilon} \bar{f}_{\varepsilon}(r) \leq -c_1 f_1(r) \,.
\end{equation}
Indeed, all the calculations in the proof of Lemma 3.6 in \cite{majka} are expressed in terms of a function $h$, which can be modified if necessary. It is enough to ensure that we choose $R_0$ such that $h^{-}(r) = 0$ for $r \geq R_0$ and then $R_1$ such that $(-r \kappa(r) + 2M)/2 \leq -C_{\varepsilon} r / (R_1 - R_0)R_1$ for $r \geq R_1$, which obviously holds for the choice of $h$, $R_0$ and $R_1$ presented above. Moreover, we obviously have
\begin{equation*}
 -f_1'(r)\kappa(r)r + C_{\varepsilon} \bar{f}_{\varepsilon}(r) \leq  -f_1'(r)\kappa(r)r + 2f_1'(r)M + C_{\varepsilon} \bar{f}_{\varepsilon}(r) \,,
\end{equation*}
and thus if we have (\ref{functionalInequalityWithM}) with some $M > 0$ for $r > \delta$, then (\ref{functionalInequality}) is still valid for $r > \delta$. The reason we introduce the constant $M$ is that it is needed to show that 
\begin{equation*}
 \sup_{x \in (r, r+\varepsilon)} f_1''(x) \leq -\frac{c_1}{C_{\varepsilon}}f_1(r) + f_1'(r)\frac{r\kappa(r)}{C_{\varepsilon}} 
\end{equation*}
holds for all $r \leq \delta$, which, combined with (\ref{functionalInequality}) for $r > \delta$, will give us (\ref{newFunctionalInequality}). Hence we need to show that for any $s \in (r, r+ \varepsilon)$ we have
\begin{equation*}
 f_1''(s) \leq -\frac{c_1}{C_{\varepsilon}}f_1(r) + f_1'(r)\frac{r\kappa(r)}{C_{\varepsilon}} \,.
\end{equation*}
First let us calculate (recall that we can choose $R_1$ large enough so that $s < \delta + \varepsilon < R_1$)
\begin{equation*}
\begin{split}
 f_1''(s) &= \phi(s) \left( - \frac{c_1}{C_{\varepsilon}} \frac{\Phi(s+\varepsilon)}{\phi(s)} \right) + \left( - \frac{\bar{h}(r)}{C_{\varepsilon}}\phi(s) g(s) \right) \\
 &= - \frac{c_1}{C_{\varepsilon}} \Phi(s + \varepsilon) - \frac{\bar{h}(r)}{C_{\varepsilon}} f_1'(s) \,.
 \end{split}
\end{equation*}
Observe now that for $r \leq s$ we have $f_1(r) \leq f_1(s) \leq \Phi(s) \leq \Phi(s + \varepsilon)$ and thus
\begin{equation*}
 f_1''(s) \leq - \frac{c_1}{C_{\varepsilon}} f_1(r) - \frac{\bar{h}(r)}{C_{\varepsilon}} f_1'(s) \,.
\end{equation*}
Therefore it remains to be shown that
\begin{equation*}
 - \frac{\bar{h}(r)}{C_{\varepsilon}} f_1'(s) \leq f_1'(r)\frac{r\kappa(r)}{C_{\varepsilon}} \,.
\end{equation*}
Actually, we will just show that
\begin{equation}\label{goal}
 \frac{1}{C_{\varepsilon}}f_1'(s)h(s)  \leq f_1'(r)\frac{r\kappa(r)}{C_{\varepsilon}} \,.
\end{equation}
Then, since $-h^{-} \leq h$ and $s \in (r, r+ \varepsilon)$, we will get
\begin{equation*}
\begin{split}
- \frac{1}{C_{\varepsilon}} f_1'(s) \sup_{t \in (r,r+\varepsilon)}h^{-}(t) &= \frac{1}{C_{\varepsilon}} f_1'(s) \inf_{t \in (r,r+\varepsilon)} (-h^{-}(t)) \\
&\leq \frac{1}{C_{\varepsilon}} f_1'(s) \inf_{t \in (r,r+\varepsilon)} h(t) \\
&\leq f_1'(r)\frac{r\kappa(r)}{C_{\varepsilon}} \,.
\end{split}
\end{equation*}
In order to show (\ref{goal}) we observe that straight from the definition of $h$ we have
\begin{equation*}
 \frac{1}{C_{\varepsilon}}f_1'(s)h(s) =  \frac{1}{C_{\varepsilon}}f_1'(s)(s\kappa(s) - 2M)
\end{equation*}
and then we calculate
\begin{equation*}
 \begin{split}
  f_1'(s)(s\kappa(s) - 2M) &= f_1'(r)r\kappa(r) - f_1'(r)r\kappa(r)\\
  &+f_1'(s)r\kappa(r) -f_1'(s)r\kappa(r) \\
  &+f_1'(s)s\kappa(s) -2M f_1'(s) \\
  &\leq f_1'(r)r\kappa(r) + r\kappa(r)(f_1'(s)-f_1'(r))\\
  &+ f_1'(s)(s\kappa(s) - r\kappa(r)) - 2M f_1'(s) \,.
  \end{split}
\end{equation*}
Now it is enough to show that it is possible to choose $\varepsilon$, $\delta$ and $M$ in such a way that the sum of the last three terms is bounded by some non-positive quantity. 
Since we have Assumption \ref{AssumptionD2}, for any $\lambda > 0$ there exists some $K(\lambda) > 0$ such that for all $|r| < \lambda$ we have $|r \kappa(r)| \leq K(\lambda)$. Since $s < r + \varepsilon \leq \delta + \varepsilon$, we obtain
\begin{equation*}
  s\kappa(s) - r\kappa(r) \leq 2 K(\delta + \varepsilon) \,.
\end{equation*}
 We also know that $f_1'$ is non-increasing and thus $f_1'(s) \leq f_1'(r)$, but the sign of $r\kappa(r)$ is unknown so we cannot just bound $r\kappa(r)(f_1'(s)-f_1'(r))$ by zero. We will deal with this term in a more complicated way. We have
\begin{equation*}
\begin{split}
 f_1'(s) - f_1'(r) &= \phi(s)g(s) - \phi(s)g(r) + \phi(s)g(r) - \phi(r)g(r) \\
 &= \phi(s)(g(s) - g(r)) + (\phi(s) - \phi(r))g(r) \,.
 \end{split}
\end{equation*}
We also have
\begin{equation*}
 |\phi(s)(g(s) - g(r))| \leq 2 \phi(s) \leq 4 f_1'(s) \,,
\end{equation*}
since $1/2 \leq g \leq 1$. Furthermore
\begin{equation*}
 \begin{split}
  |(\phi(s) - \phi(r))g(r)| &= |\phi(s)(1 - \phi(s)^{-1}\phi(r))g(r)|\\
  &= \left|\phi(s)\left(1 - \exp \left( \int_r^s \frac{\bar{h}(t)}{C_{\varepsilon}}dt \right) \right)g(r)\right|\\
  &\leq 2 f_1'(s) \int_r^s \frac{\bar{h}(t)}{C_{\varepsilon}}dt \exp \left( \int_r^s \frac{\bar{h}(t)}{C_{\varepsilon}}dt \right)\\
  &\leq  \frac{2 f_1'(s) \varepsilon}{C_{\varepsilon}}(2M + K(\delta + 2\varepsilon)) \exp\left(\frac{\varepsilon}{C_{\varepsilon}}(2M + K(\delta + 2\varepsilon))\right) ,
 \end{split}
\end{equation*}
where in the first inequality we have used the fact that $|1-e^x|\leq |xe^x|$ for all $x \geq 0$ and that $g \leq 1$ and $\phi(s) \leq 2f_1'(s)$. In the second inequality we used 
\begin{equation*}
 \int_r^s \frac{\bar{h}(t)}{C_{\varepsilon}}dt \leq \frac{\varepsilon}{C_{\varepsilon}}(2M + K(\delta + 2\varepsilon)) \,,
\end{equation*}
which holds since $|s - r| < \varepsilon$. Thus if we find $\delta$, $\varepsilon$ and $M$ such that
\begin{equation*}
 K(\delta)\Bigg(4+ 2\frac{\varepsilon}{C_{\varepsilon}}(2M + K(\delta + 2\varepsilon)) \exp\left(\frac{\varepsilon}{C_{\varepsilon}}(2M + K(\delta + 2\varepsilon))\right)\Bigg) 
 + 2K(\delta + \varepsilon) \leq 2M \,,
\end{equation*}
then (\ref{goal}) holds and we prove our statement. This is indeed possible since we assume that $\varepsilon / C_{\varepsilon}$ is bounded in a neighbourhood of zero.
\end{proof}
\end{theorem}

\section{Coupling of jump diffusions}\label{sectionCouplingJumpDiffusion}

Here we study jump diffusions of more general form (\ref{twoNoisesSDE}) and we prove Theorem \ref{mainCouplingTheorem}. In order to do this, we first recall results obtained by Eberle in \cite{eberle} for diffusions of the form
\begin{equation*}
 dX_t = b(X_t)dt + \sigma_1 dB_t \,,
\end{equation*}
where $\sigma_1$ is a constant non-degenerate $d \times d$ matrix and $(B_t)_{t \geq 0}$ is a $d$-dimensional Brownian motion. Eberle used the coupling by reflection $(X_t, Y_t)_{t \geq 0}$, defined by
\begin{equation}\label{reflectionCouplingEberle}
 dY_t = \begin{cases}
         b(Y_t)dt + \sigma_1 R_{\sigma_1}(X_t,Y_t)dB_t & \text{ for } t < T \,, \\
         dX_t & \text{ for } t \geq T \,,
        \end{cases}
\end{equation}
where $T := \inf \{ t \geq 0 : X_t = Y_t \}$ is the coupling time and
\begin{equation}\label{reflectionOperatorEberle}
 R_{\sigma_1}(X_t,Y_t) := I - 2 e_t e_t^T
\end{equation}
with
\begin{equation}\label{etEberle}
 e_t := \sigma_1^{-1}(X_t - Y_t) / |\sigma_1^{-1}(X_t - Y_t)| \,.
\end{equation}
Using this coupling, Eberle constructed a concave continuous function $f$ given by
\begin{equation}\label{defFEberle}
 f(r) := \int_0^r \varphi(s) g(s) ds \,,
\end{equation}
where
\begin{align*}
 \varphi(r) &:= \exp \left( -\frac{1}{2} \int_0^r s \kappa^{-}(s) ds \right) \,, &\quad   g(r) &:= \begin{cases} 
      1 - \frac{\alpha c}{2} \int_0^r \frac{\Phi(t)}{\varphi(t)} dt \,, & r\leq R_1 \,, \\
      \frac{1}{2} \,, & r\geq R_1 \,,
   \end{cases}
\end{align*}
with $\Phi(r) := \int_0^r \varphi(s) ds$ and some constant $R_1 > 0$ defined by (9) in \cite{eberle}. Here $c > 0$ is a~constant given by
\begin{equation}\label{defcEberle}
 c = \frac{1}{\alpha} \left( \int_0^{R_1} \frac{\Phi(s)}{\varphi(s)} ds \right)^{-1}
\end{equation}
where $\alpha := \sup \{ |\sigma_1^{-1}z|^2 : z \in \mathbb{R}^d \text{ with } \| z \| = 1 \}$ (cf. the formula (12) in \cite{eberle}) and $\kappa$ is defined by
\begin{equation}\label{kappaEberle}
 \kappa(r) = \inf \left\{ -\frac{|\sigma_1^{-1}(x-y)|^2}{|x-y|^4}\langle b(x) - b(y), x - y \rangle : x, y \in \mathbb{R}^d \text{ s.t. } |x-y| = r \right\} \,.
\end{equation}
In other words, $\kappa$ is the largest quantity satisfying
\begin{equation}\label{kappaEberleInequality}
 \langle b(x) - b(y), x - y \rangle \leq -\kappa(|x-y|) |x-y|^4 / |\sigma_1^{-1}(x-y)|^2
\end{equation}
for all $x$, $y \in \mathbb{R}^d$, although for our purposes we can consider any continuous function $\kappa$ such that (\ref{kappaEberleInequality}) holds. Then it is possible to prove that
\begin{equation}\label{functionalInequalityEberle}
 2f''(r) - r\kappa(r) f'(r) \leq -c \alpha f(r) \text{ for all } r > 0 \,.
\end{equation}
Note that our definition of $\kappa$ differs from the one in \cite{eberle} by a factor $2$ to make the notation more consistent with our results for the pure jump noise case presented in the previous Section (cf. formulas in Section 2.1 in \cite{eberle}). By the methods explained in the proof of Theorem 1 in \cite{eberle} (see also Corollary 2 therein) we get
\begin{equation*}
 \mathbb{E}f(|X_t(x) - Y_t(y)|) \leq e^{-ct} f(|x-y|)
\end{equation*}
and, by the choice of $f$ (which is comparable with the identity function, since it is extended in an affine way from $R_1 > 0$), we also get
\begin{equation*}
 \mathbb{E}|X_t(x) - Y_t(y)| \leq C e^{-ct} |x-y|
\end{equation*}
with a constant $C > 0$ defined by (cf. (14) and (8) in \cite{eberle})
\begin{equation}\label{defCEberle}
 C:= 2 \varphi(R_0)^{-1} \text{, where } R_0 := \inf \{ R \geq 0 : \forall r \geq R \; \kappa (r) \geq 0 \} \,.
\end{equation}

Now we will explain how to combine the results from \cite{eberle} and \cite{majka} to get analogous results for equations involving both the Gaussian and the Poissonian noise. The general idea is, similarly to \cite{eberle} and \cite{majka}, to use an appropriate coupling $(X_t,Y_t)_{t \geq 0}$, to write an SDE for the difference process $Z_t = X_t - Y_t$, to use the It\^{o} formula to evaluate $df(|Z_t|)$ and then to choose $f$ in such a way that $df(|Z_t|) \leq dM_t - \widetilde{c}f(|Z_t|)dt$ for some constant $\widetilde{c} > 0$, where $(M_t)_{t \geq 0}$ is a local martingale.

\begin{proofMainCoupling}

We consider an equation of the form
\begin{equation*}
 dX_t = b(X_t)dt + \sigma_1 dB^1_t + \sigma(X_t) dB^2_t + dL_t + \int_{U} g(X_{t-},u) \widetilde{N}(dt,du) \,,
\end{equation*}
where $\sigma_1 > 0$ is a constant and all the other coefficients and the sources of noise are like in the formulation of Theorem \ref{mainCouplingTheorem} (in particular, here we denote the underlying Poisson random measure of $(L_t)_{t \geq 0}$ by $N^L$ and its associated L\'{e}vy measure by $\nu^L$). Restricting ourselves to a real constant in front of $(B^1_t)_{t \geq 0}$ instead of a matrix helps us to slightly reduce the notational complexity and seems in fact quite natural at least for the equations for which Lemma \ref{diffusionSplittingLemma} applies. Recall that $\kappa$ is such that for all $x$, $y \in \mathbb{R}^d$ we have
\begin{equation}\label{kappaSect4}
 \langle b(x) - b(y) , x - y \rangle + \frac{1}{2}\int_U |g(x,u) - g(y,u)|^2 \nu(du) + \| \sigma(x) - \sigma(y) \|_{HS}^2 \leq -\kappa(|x-y|)|x-y|^2
\end{equation}
and that it satisfies Assumption \ref{AssumptionD1}.
Now we will apply the mirror coupling from \cite{majka} to $(L_t)_{t \geq 0}$, by using the ``mirror operator'' $M(\cdot, \cdot)$, i.e., recalling the notation used in the equations (\ref{eqX}) and (\ref{optimalY}), we define
\begin{equation*}
 \begin{split}
  M(X_{t-},Y_{t-})L_t &:= \int_{\{ |v| > m \} \times [0,1]} v N^L(dt,dv,du) \\
  &+ \int_{\{ |v| \leq m \} \times [0,1]} (X_{t-} - Y_{t-} + v)\mathbf{1}_{\{ u < \rho(v, Z_{t-}) \}} \widetilde{N^L}(dt,dv,du) \\
  &+ \int_{\{ |v| \leq m \} \times [0,1]} R(X_{t-},Y_{t-})v\mathbf{1}_{\{ u \geq \rho(v, Z_{t-}) \}} \widetilde{N^L}(dt,dv,du) \,,
 \end{split}
\end{equation*}
with the reflection operator $R$ defined by (\ref{reflectionOperatorLevy}). We will also use the reflection coupling (\ref{reflectionCouplingEberle}) from \cite{eberle}, with the reflection operator $R_{\sigma_1}$ defined by (\ref{reflectionOperatorEberle}) and apply it to $(B^1_t)_{t \geq 0}$. Note that if the coefficient near the Brownian motion is just a positive constant and not a matrix, the formulas from \cite{eberle} become a bit simpler, in particular the unit vector $e_t$ defined by (\ref{etEberle}) becomes just $(X_t - Y_t) / |X_t - Y_t|$. Thus the two reflection operators we defined coincide and we can keep denoting them both by $R$. Moreover, we apply the synchronous coupling to the other two noises and hence we have
\begin{equation*}
 dY_t = b(Y_t)dt + \sigma_1 R(X_t,Y_t) dB^1_t + \sigma(Y_t) dB^2_t + M(X_{t-},Y_{t-}) dL_t + \int_{U} g(Y_{t-},u) \widetilde{N}(dt,du) \,.
\end{equation*}
Since all the sources of noise are independent, it is easy to see that $(X_t,Y_t)_{t \geq 0}$ is indeed a~coupling (this follows from the fact that $R$ applied to $(B_t^1)_{t \geq 0}$ gives a Brownian motion and $M$ applied to $(L_t)_{t \geq 0}$ gives the same L\'{e}vy process, whereas the solution to the equation above is unique in law).
We can now write the equation for $Z_t := X_t - Y_t$ as
\begin{equation*}
\begin{split}
 dZ_t &= (b(X_t) - b(Y_t))dt + 2\sigma_1 e_t e_t^T dB^1_t + (\sigma(X_t) - \sigma(Y_t))dB^2_t \\
 &+ (I - M(X_{t-},Y_{t-})) dL_t +  \int_{U} (g(X_{t-},u) - g(Y_{t-},u)) \widetilde{N}(dt,du) \,,
 \end{split}
\end{equation*}
where we evaluated $\sigma_1(I - R(X_{t-},Y_{t-}))$ as $2\sigma_1 e_t e_t^T$ and we will later use the fact that $d\widetilde{W}_t := e_t^T dB^1_t$ is a one-dimensional Brownian motion in order to simplify our calculations. We apply the It\^{o} formula to get
\begin{equation}\label{ItoSum}
 df(|Z_t|) = \sum_{j=1}^9 I_j \,,
\end{equation}
where
\begin{align*}
 I_1 &:= f'(|Z_t|)\frac{1}{|Z_t|} \langle b(X_t) - b(Y_t) , Z_t \rangle dt \,, &\,  I_3 &:= f'(|Z_t|)\frac{1}{|Z_t|} \langle  Z_t , (\sigma(X_t) - \sigma(Y_t))dB^2_t \rangle \,, \\
 I_2 &:= 2f'(|Z_t|)\frac{1}{|Z_t|} \langle Z_t , \sigma_1 e_t e_t^T dB^1_t \rangle \,,  &\, I_4 &:= f'(|Z_{t-}|)\frac{1}{|Z_{t-}|} \langle Z_{t-} , (I - M(X_{t-},Y_{t-}))  \rangle dL_t
\end{align*}
and
\begin{equation*}
 I_5 := f'(|Z_{t-}|)\frac{1}{|Z_{t-}|} \int_U \langle g(X_{t-},u) - g(Y_{t-},u) , Z_{t-} \rangle \widetilde{N}(dt,du)
\end{equation*}
constitute the drift and the local martingale terms, while
\begin{equation*}
 I_6 := \frac{1}{2} \sigma_1^2 \sum_{i, j = 1}^{d} \Bigg[ f''(|Z_{t-}|) \frac{Z_{t-}^i Z_{t-}^j}{|Z_{t-}|^2} + f'(|Z_{t-}|)(\delta_{ij} \frac{1}{|Z_{t-}|} - \frac{Z_{t-}^i Z_{t-}^j}{|Z_{t-}|^3} ) \Bigg] 4 \frac{Z_{t-}^i Z_{t-}^j}{|Z_{t-}|^2} dt
\end{equation*}
and
\begin{equation*}
 \begin{split}
  I_7 &:= \sum_{i, j = 1}^{d} \Bigg[ f''(|Z_{t-}|) \frac{Z_{t-}^i Z_{t-}^j}{|Z_{t-}|^2} + f'(|Z_{t-}|)(\delta_{ij} \frac{1}{|Z_{t-}|} - \frac{Z_{t-}^i Z_{t-}^j}{|Z_{t-}|^3} ) \Bigg] \\
  &\cdot \Bigg[ \sum_{k = 1}^{m} (\sigma_{ik}(X_{t-}) - \sigma_{ik}(Y_{t-}))(\sigma_{jk}(X_{t-}) - \sigma_{jk}(Y_{t-})) \Bigg] dt
 \end{split}
\end{equation*}
come from the quadratic variation of the Brownian noises, whereas
\begin{equation*}
\begin{split}
 I_8 &:= \int_U \Bigg[ f(|Z_{t-} + (I - M(X_{t-},Y_{t-}))v|) - f(|Z_{t-}|) \\ 
 &- \langle (I - M(X_{t-},Y_{t-}))v , \nabla f(|Z_{t-}|) \rangle \Bigg] N^L(dt,dv)
 \end{split}
\end{equation*}
and
\begin{equation}\label{I9}
\begin{split}
 I_9 &:= \int_U \Bigg[ f(|Z_{t-} + g(X_{t-},u) - g(Y_{t-},u)|) - f(|Z_{t-}|) \\
 &- \langle g(X_{t-},u) - g(Y_{t-},u) , \nabla f(|Z_{t-}|) \rangle \Bigg] N(dt,du)
 \end{split}
\end{equation}
are the jump components.

Now we proceed similarly to \cite{eberle} and \cite{majka}. Since we want to obtain an estimate of the form $df(|Z_t|) \leq dM_t - \widetilde{c}f(|Z_t|)dt$ and we assume that the function $f$ is concave, we should use its second derivative to obtain a negative term on the right hand side of (\ref{ItoSum}). In order to do this, we can use the additive Brownian noise $(B^1_t)_{t \geq 0}$ to get a negative term from $I_6$ (it is easy to see that it reduces to $2 \sigma_1^2 f''(|Z_t|)$) and then use the function $f$ from \cite{eberle} given by (\ref{defFEberle}), aiming to obtain an inequality like (\ref{functionalInequalityEberle}) (then we can just use the synchronous coupling for $(L_t)_{t \geq 0}$ and the terms $I_4$ and $I_8$ disappear). Alternatively, we can use the additive jump noise $(L_t)_{t \geq 0}$ to get a negative term from $I_8$. As we already mentioned in Section \ref{sectionCouplingPureJump} under the formula (\ref{ItoOptimal}), the integral $I_8$ reduces to the left hand side of (\ref{ItoJumpSum}), see Section 3 in \cite{majka} for details. Then we can use the function $f_1$ from \cite{majka}, aiming to obtain an inequality like (\ref{newFunctionalInequality}) (then we use the synchronous coupling for $(B^1_t)_{t \geq 0}$ and the terms $I_2$ and $I_6$ disappear). In either case, $I_3$ and $I_5$ can be controlled via $\kappa$, since the coefficients $\sigma$ and $g$ are included in its definition. If we are only interested in finding any constant $\widetilde{c} > 0$ such that (\ref{initialPerturbationInequality}) holds, then it is sufficient to use one of the two additive noises and to apply the synchronous coupling to the other (if both noises are present it is recommendable to use $(B^1_t)_{t \geq 0}$ since the formulas in \cite{eberle} are simpler than the ones in \cite{majka}). If we are interested in finding the best (largest) possible constant $\widetilde{c}$, then we can use both noises, but then we would also need to redefine the function $f$ and this would be technically quite sophisticated (whereas by using only one noise we can essentially just use the formulas that are already available in either \cite{eberle} or \cite{majka}).

We should still explain how to control $I_7$ and $I_9$. We can control $I_9$ following the ideas from \cite{luowang} and controlling $I_7$ is also quite straightforward.

First observe that $\nabla f(|Z_{t-}|) = f'(|Z_{t-}|) \frac{1}{|Z_{t-}|} Z_{t-}$ and, following Section 5.2 in \cite{luowang}, note that since $f$ is concave and differentiable, we have
\begin{equation*}
 f(a) - f(b) \leq f'(b)(a-b)
\end{equation*}
for any $a$, $b > 0$. Thus
\begin{equation*}
\begin{split}
 &f(|Z_{t-} + g(X_{t-},u) - g(Y_{t-},u)|) - f(|Z_{t-}|) - f'(|Z_{t-}|) \frac{1}{|Z_{t-}|} \langle g(X_{t-},u) - g(Y_{t-},u) , Z_{t-} \rangle \\
 &\leq f'(|Z_{t-}|) \left( |Z_{t-} + g(X_{t-},u) - g(Y_{t-},u)| - |Z_{t-}| - \frac{1}{|Z_{t-}|} \langle g(X_{t-},u) - g(Y_{t-},u) , Z_{t-} \rangle \right) \,.
\end{split}
\end{equation*}
Next we will need the inequality
\begin{equation*}
|x+y| - |x| - \frac{1}{|x|}\langle y, x \rangle \leq \frac{1}{2|x|}|y|^2 \,,
\end{equation*}
which holds for any $x$, $y \in \mathbb{R}^d$ since
\begin{equation*}
 |x| |x+y| \leq \frac{1}{2}(|x|^2 + |x+y|^2) = \frac{1}{2}(|x|^2 + |x|^2 + 2\langle x , y \rangle + |y|^2) \,.
\end{equation*}
Hence we obtain 
\begin{equation*}
 I_9 \leq f'(|Z_{t-}|) \frac{1}{2|Z_{t-}|} \int_U |g(X_{t-},u) - g(Y_{t-},u)|^2 N(dt,du) \,.
\end{equation*}

On the other hand, if we denote by $\sigma^k$ the $k$-th column of the matrix $\sigma$, then
\begin{equation*}
\begin{split}
 I_7 &= \sum_{k=1}^{m} f''(|Z_{t-}|) \frac{|\langle Z_t , \sigma^k(X_t) - \sigma^k(Y_t) \rangle |^2}{|Z_{t-}|^2} +  \sum_{k=1}^{m} \sum_{i=1}^{d} f'(|Z_{t-}|) \frac{1}{|Z_{t-}|}(\sigma_{ik}(X_{t-}) - \sigma_{ik}(Y_{t-}))^2 \\
 &- \sum_{k=1}^{m} f'(|Z_{t-}|) \frac{|\langle Z_t , \sigma^k(X_t) - \sigma^k(Y_t) \rangle |^2}{|Z_{t-}|^3}  \leq f'(|Z_{t-}|) \frac{1}{|Z_{t-}|} \| \sigma(X_{t-}) - \sigma(Y_{t-}) \|_{HS}^2 \,.
 \end{split}
\end{equation*}
Hence we get a bound on $df(|Z_t|)$, which allows us to bound $\mathbb{E}f(|Z_t|) - \mathbb{E}f(|Z_s|)$ for any $0 \leq s < t$. Using a localization argument with a sequence of stopping times $(\tau_n)_{n=1}^{\infty}$ like in (\ref{stoppingTimes}), we can get rid of the expectations of the local martingale terms. Then we can use the inequality (\ref{kappaSect4}) multiplied by $f'(|x-y|)\frac{1}{|x-y|}$ to see that, if we are using the additive L\'{e}vy noise $(L_t)_{t \geq 0}$ to get our bounds, then after handling $I_8$ like in (\ref{ItoJumpSum}) and using the estimates (\ref{majkaBound1}) and (\ref{majkaBound2}), we need to choose a function $f_1$~such that
\begin{equation*}
-f_1'(r)\kappa(r)r + C_{\varepsilon} \hat{f}_{\varepsilon}(r) \leq -c_1 f_1(r)
\end{equation*}
and this is exactly (\ref{newFunctionalInequality}), so we can handle further calculations like in the proof of Theorem \ref{jumpCouplingMainTheorem}. Alternatively, if we are using the additive Gaussian noise $(B_t^1)_{t \geq 0}$, we can modify the definition of $\kappa$ to include the $\sigma_1$ factor (cf. (\ref{kappaEberle})) and then we need to choose a function $f$ such that
\begin{equation*}
 2 f''(r) - r\kappa(r) f'(r) \leq -\frac{c}{\sigma_1^2} f(r) \,,
\end{equation*}
hence $1 / \sigma_1^2$ plays the role of $\alpha$ in the calculations in \cite{eberle} (cf. (\ref{functionalInequalityEberle}) earlier in this Section and for the details see the proof of Theorem 1 in \cite{eberle}, specifically the formula (63), while remembering about the change of the factor $2$ in our definition of $\kappa$ compared to the one in \cite{eberle}).

Either way we obtain some constant $\widetilde{c} > 0$ and a function $\widetilde{f}$ such that 
\begin{equation}\label{preGronwall}
 \mathbb{E}\widetilde{f}(|Z_{t \wedge \tau_n}|) - \mathbb{E}\widetilde{f}(|Z_{s \wedge \tau_n}|) \leq -\widetilde{c} \int_s^t \mathbb{E}\widetilde{f}(|Z_{r \wedge \tau_n}|)dr
\end{equation}
holds for any $0 \leq s < t$. Here $\widetilde{c}$ and $\widetilde{f}$ are equal either to $c$ and $f$ defined by (\ref{defcEberle}) and (\ref{defFEberle}) or $c_1$ and $f_1$ defined by (\ref{defc1}) and (\ref{defFMajka}), respectively, depending on whether we used $(B_t^1)_{t \geq 0}$ or $(L_t)_{t \geq 0}$ in the step above. Thus we can use the differential version of the Gronwall inequality to get
\begin{equation*}
 \mathbb{E}\widetilde{f}(|Z_{t \wedge \tau_n}|) \leq \mathbb{E}\widetilde{f}(|Z_0|) e^{-\widetilde{c}t} \text{ for any } t > 0 \,,
\end{equation*}
and after using the Fatou lemma, the fact that $\tau_n \to T$ and that $Z_t = 0$ for $t \geq T$, we get 
\begin{equation*}
 \mathbb{E}\widetilde{f}(|Z_t|) \leq \mathbb{E}\widetilde{f}(|Z_0|) e^{-\widetilde{c}t} \text{ for any } t > 0 \,.
\end{equation*}
Since we can compare our function $\widetilde{f}$ with the identity function from both sides, this finishes the proof. Note that in the last step one has to be careful and use the differential version of the Gronwall formula, since the integral version does not work when the term on the right hand side is negative (cf. Remark 2.3 in \cite{shao}). 

Note also that if we are only dealing with the Gaussian noise, then we can reason like in \cite{eberle}, i.e., having proved that
\begin{equation*}
 df(|Z_t|) \leq dM_t - cf(|Z_t|)dt
\end{equation*}
(by choosing an appropriate function $f$) for some local martingale $(M_t)_{t \geq 0}$, we can see that this implies $d(e^{ct} f(|Z_t|)) \leq dM_t$, so by using a localization argument we can directly get $\mathbb{E}[e^{ct}f(|Z_t|)] \leq \mathbb{E}f(|Z_0|)$ without using the Gronwall inequality. However, in the jump case this is not possible, since we have to first take the expectation in order to deal with $I_8$ and $I_9$ by transforming the stochastic integrals with respect to $N^L$ and $N$ into deterministic integrals with respect to $\nu^L$ and $\nu$, respectively. Only then can we use the definition of $\kappa$ via (\ref{kappaSect4}) to find an appropriate function $\widetilde{f}$ such that (\ref{preGronwall}) holds.

\qed

\end{proofMainCoupling}

We will now show how, starting from an equation of the form (\ref{generalSDEsect2}) with one multiplicative Gaussian noise, we can obtain an SDE of the form (\ref{twoNoisesSDE}) with two independent Gaussian noises, one of which is still multiplicative, but the other additive (and the additive one has just a real constant as a coefficient, cf. the comments in the proof of Theorem \ref{mainCouplingTheorem} earlier in this section).

\begin{lemma}\label{diffusionSplittingLemma}
If $(X_t)_{t \geq 0}$ is the unique strong solution to the SDE
\begin{equation}\label{multiplicativeDiffusion}
 dX_t = b(X_t)dt + \sigma(X_t)dB_t \,,
\end{equation}
where $(B_t)_{t \geq 0}$ is a Brownian motion and $\sigma \sigma^T$ is uniformly positive definite, then $(X_t)_{t \geq 0}$ can also be obtained as a solution to
\begin{equation}\label{splitDiffusion}
 dX_t = b(X_t)dt + CdB_t^1 + \widetilde{\sigma}(X_t)dB_t^2
\end{equation}
with two independent Brownian motions $(B_t^1)_{t \geq 0}$ and $(B_t^2)_{t \geq 0}$, some constant $C > 0$ and a diffusion coefficient $\widetilde{\sigma}$ such that if $\| \sigma(x) \|_{HS} \leq M$ for all $x \in \mathbb{R}^d$ with some constant $M > 0$, then
\begin{equation}\label{LipschitzOfModifiedSigma}
 \| \widetilde{\sigma}(x) - \widetilde{\sigma}(y)\|_{HS} \leq \frac{M}{\sqrt{\lambda^2 - C^2}} \| \sigma(x) - \sigma(y) \|_{HS} \,,
\end{equation}
where the constants $\lambda > C > 0$ are as indicated in the proof.

\begin{proof}
 Observe that if the diffusion coefficient $\sigma$ is such that $\sigma \sigma^T$ is uniformly positive definite, i.e., there exists $\lambda > 0$ such that for any $x$, $h \in \mathbb{R}^d$ we have
\begin{equation*}
 \langle \sigma(x)\sigma(x)^T h , h \rangle \geq \lambda^2 |h|^2 \,,
\end{equation*}
then $\sigma(x)\sigma(x)^T - \lambda^2 I$ is nonnegative definite for any $x \in \mathbb{R}^d$ and thus we can consider
\begin{equation*}
 a(x) := \sqrt{\sigma(x)\sigma(x)^T - \lambda^2 I} \,,
\end{equation*}
which is the unique (symmetric) nonnegative definite matrix such that $a(x)a(x)^T = \sigma(x)\sigma(x)^T - \lambda^2 I$. Note that if we now define $\widetilde{\sigma}$ as
\begin{equation*}
 \widetilde{\sigma}(x) := \sqrt{\sigma(x)\sigma(x)^T - C^2 I} 
\end{equation*}
for some constant $0 < C^2 < \lambda^2$, we can get
\begin{equation*}
 \langle \widetilde{\sigma}(x)^2 h , h \rangle = \langle \sigma(x)\sigma(x)^T h , h \rangle - C^2 \langle h , h \rangle \geq (\lambda^2 - C^2)|h|^2 \,, 
\end{equation*}
and thus we can assume that $\widetilde{\sigma}(x)$ is also uniformly positive definite. Therefore Lemma 3.3 in \cite{priola} applies (our $\widetilde{\sigma}$ corresponds to $\sigma$ in \cite{priola} and our $\sigma \sigma^T$ corresponds to $q$ therein). Thus we get
\begin{equation*}
 \| \widetilde{\sigma}(x) - \widetilde{\sigma}(y) \|_{HS} \leq \frac{1}{2\sqrt{\lambda^2 - C^2}} \| \sigma(x) \sigma(x)^T - \sigma(y) \sigma(y)^T \|_{HS} \,.
\end{equation*}
(all eigenvalues of $\sigma(x)\sigma(x)^T - C^2 I$ are not less than $\lambda^2 - C^2$, which is the condition that needs to be checked in the proof of Lemma 3.3 in \cite{priola}). This shows that whenever $\sigma \sigma^T$ is Lipschitz with a constant $L$, the function $\widetilde{\sigma}$ is Lipschitz with $L / 2\sqrt{\lambda^2 - C^2}$. In particular, if $\sigma$ is Lipschitz with a constant $L$ and bounded with a constant $M$, then $\sigma \sigma^T$ is Lipschitz with the constant $2LM$ and thus $\widetilde{\sigma}$ is Lipschitz with  $LM/\sqrt{\lambda^2 - C^2} $. Hence we prove (\ref{LipschitzOfModifiedSigma}).

Now assume that $(X_t)_{t \geq 0}$ is a solution to (\ref{splitDiffusion}) and consider the process
\begin{equation*}
 A_t := C B^1_t + \int_0^t \sqrt{\sigma(X_s)\sigma(X_s)^T - C^2 I} dB^2_s = X_t - X_0 - \int_0^t b(X_s) ds \,.
\end{equation*}
We can easily calculate
\begin{equation}\label{quadraticVariation}
 [A^i,A^j]_t = \int_0^t \left( \sigma \sigma^T \right)_{ij} (X_s) ds \,.
\end{equation}
Hence, if we write
\begin{equation*}
 dX_t = dA_t + b(X_t)dt = \sigma(X_t) d\widetilde{B}_t + b(X_t)dt \,,
\end{equation*}
where $d\widetilde{B}_t =  \sigma^{-1}(X_t)dA_t$, then using (\ref{quadraticVariation}) and the L\'{e}vy characterization theorem, we infer that $(\widetilde{B}_t)_{t \geq 0}$ is a Brownian motion. Thus $(X_t)_{t \geq 0}$ is a solution to (\ref{multiplicativeDiffusion}).
\end{proof}
\end{lemma}

The proof of (\ref{LipschitzOfModifiedSigma}) is based on the reasoning in \cite{priola}, Section 3 (the matrix $q(x)$ used there is our $\sigma(x)\sigma(x)^T$; the difference in notation follows from the fact that the starting point for studying diffusions in \cite{priola} is the generator and not the SDE). Due to (\ref{LipschitzOfModifiedSigma}) we see that if the coefficients in (\ref{multiplicativeDiffusion}) satisfy Assumption \ref{AssumptionD1}, then the coefficients in the modified equation (\ref{splitDiffusion}) also do (after a suitable change in the definition of $\kappa$). More generally, Lemma \ref{diffusionSplittingLemma} allows us to replace an equation of the form (\ref{generalSDEsect2}) with
\begin{equation*}
 dX_t = b(X_t)dt + C dB^1_t + \sqrt{\sigma(X_t)\sigma(X_t)^T - C^2 I} dB^2_t + \int_{U} g(X_{t-},u) \widetilde{N}(dt,du) \,,
\end{equation*}
as long as $\sigma \sigma^T$ is uniformly positive definite.

\section{Bounds on Malliavin derivatives}\label{sectionMalliavin}

\subsection{Brownian case}\label{subsectionBrownianMalliavin}

In this section we first prove Theorem \ref{theoremDriftPerturbation} and then we show how to obtain bounds on Malliavin derivatives using the inequality (\ref{driftChangeBound}). As a consequence we prove Theorem \ref{generalMalliavinResult} and Corollary \ref{corollaryMalliavinBounds}. We begin with proving the following crucial result.

\begin{lemma}\label{lemmaForMalliavin}
 Let $(X_t)_{t \geq 0}$ be a $d$-dimensional jump diffusion process given by
 \begin{equation}\label{MalliavinEstimateSDE}
  dX_t = b(X_t)dt + \sigma dB_t + \int_U g(X_{t-},u) \widetilde{N}(dt,du) \,,
 \end{equation}
where $\sigma$ is a $d \times d$ matrix with $\det \sigma > 0$ and $(B_t)_{t \geq 0}$ is a $d$-dimensional Brownian motion, whereas $b$~and $g$ satisfy Assumption \ref{AssumptionD1} and Assumption \ref{AssumptionD2} and $g$ is Lipschitz. Let $h_t$ be an adapted $d$-dimensional process and consider a jump diffusion $(\widetilde{X}_t)_{t \geq 0}$ with the drift perturbed by $h_t$, i.e.,
\begin{equation}\label{Xtilde}
 d\widetilde{X}_t = b(\widetilde{X}_t)dt + h_t dt + \sigma dB_t + \int_U g(\widetilde{X}_{t-},u) \widetilde{N}(dt,du)\,.
\end{equation}
Then there exists a $d$-dimensional process $(Y_t)_{t \geq 0}$ such that $(X_t, Y_t)_{t \geq 0}$ is a coupling and we have
\begin{equation}\label{XYestimate}
 \mathbb{E} |\widetilde{X}_t - Y_t| \leq C \mathbb{E} \int_0^t e^{c(s-t)} |h_s| ds 
\end{equation}
for some constant $C > 0$.
\begin{proof}
 The arguments we use here are based on ideas from Sections 6 and 7 in \cite{eberle}, where interacting diffusions (without the jump noise) were studied. Here the most important part of the argument also concerns the Gaussian noise, however, we include the jump noise too in order to show how to handle the additional terms, which is important for the proof of Theorem \ref{theoremDriftPerturbation}. On the other hand, in order to slightly simplify the notation, we assume from now on that $\sigma = I$. Denote
 \begin{equation*}
  Z_t := \widetilde{X}_t - Y_t \,,
 \end{equation*}
where $(Y_t)_{t \geq 0}$ will be defined below by (\ref{MalliavinEstimateY}), and consider Lipschitz continuous functions $\lambda$, $\pi: \mathbb{R}^d \to [0,1]$ such that for some fixed $\delta > 0$ we have
\begin{equation}\label{lambdaProperties}
\begin{split}
 \lambda^2(z) + \pi^2(z) &= 1 \text{ for any } z \in \mathbb{R}^d \,, \\
 \lambda(z) &= 0 \text{ if } |z| \leq \delta / 2 \\
 \lambda(z) &= 1 \text{ if } |z| \geq \delta \,.
 \end{split}
\end{equation}
Now fix a unit vector $u \in \mathbb{R}^d$ and define $R(\widetilde{X}_t, Y_t) := I - 2e_t e_t^T$, where
\begin{equation*}
e_t :=  \begin{cases}
  \frac{Z_t}{|Z_t|} \,, & \text{ if } \widetilde{X}_t \neq Y_t \,, \\
  u \,, & \text{ if } \widetilde{X}_t = Y_t \,.
 \end{cases}
\end{equation*}
We will see from the proof that the exact value of $u$ is irrelevant. Let us notice that the equation (\ref{Xtilde}) for the process $(\widetilde{X}_t)_{t \geq 0}$ can be rewritten as
\begin{equation}\label{XtildeNew}
 d\widetilde{X}_t = b(\widetilde{X}_t)dt + h_t dt + \lambda(Z_t) dB^1_t + \pi(Z_t) dB^2_t + \int_U g(\widetilde{X}_{t-},u) \widetilde{N}(dt,du)\,,
\end{equation}
where $(B^1_t)_{t \geq 0}$ and $(B^2_t)_{t \geq 0}$ are independent Brownian motions, and define
\begin{equation}\label{MalliavinEstimateY}
 dY_t = b(Y_t)dt + \lambda(Z_t) R(\widetilde{X}_t, Y_t) dB^1_t + \pi(Z_t) dB^2_t + \int_U g(Y_{t-},u) \widetilde{N}(dt,du)\,.
\end{equation}
Using the L\'{e}vy characterization theorem and the fact that $\lambda^2 + \pi^2 = 1$, we can show that the processes defined by
\begin{equation*}
\begin{split}
 d\widetilde{B}_t &:= \lambda(Z_t) dB^1_t + \pi(Z_t) dB^2_t \,, \\
 d\bar{B}_t &:= \lambda(Z_t) R(\widetilde{X}_t, Y_t) dB^1_t + \pi(Z_t) dB^2_t 
\end{split}
\end{equation*}
are both $d$-dimensional Brownian motions and hence the process $(Y_t)_{t \geq 0}$ defined by (\ref{MalliavinEstimateY}) has the same finite dimensional distributions as $(X_t)_{t \geq 0}$ defined by (\ref{MalliavinEstimateSDE}) with $\sigma = I$, while both (\ref{Xtilde}) with $\sigma = I$ and (\ref{XtildeNew}) also define the same (in law) process, which follows from the uniqueness in law of solutions to equations of the form (\ref{MalliavinEstimateSDE}). Thus $(X_t, Y_t)_{t \geq 0}$ is a~coupling. Note that obviously in this case $(\widetilde{X}_t, Y_t)_{t \geq 0}$ is not a coupling, but we do not need this to prove (\ref{XYestimate}). Consider the equation for $Z_t = \widetilde{X}_t - Y_t$, which is given by
\begin{equation*}
 dZ_t = (b(\widetilde{X}_t) - b(Y_t))dt + h_t dt + 2 \lambda(Z_t) e_t e_t^T dB^1_t + \int_U (g(\widetilde{X}_{t-},u) - g(Y_{t-},u)) \widetilde{N}(dt,du)\,,
\end{equation*}
and observe that the process
\begin{equation*}
 d\widetilde{W}_t := e_t^T dB^1_t
\end{equation*}
is a one-dimensional Brownian motion. Now we would like to apply the It\^{o} formula to calculate $df(|Z_t|)$ for the function $f$ given by (\ref{defFEberle}), just like we did in the proof of Theorem \ref{mainCouplingTheorem}. However, the function $x \mapsto f(|x|)$ is not differentiable at zero. In the proof of Theorem \ref{mainCouplingTheorem} this was not a problem, since we started the marginal processes of our coupling at two different initial points and were only interested in the behaviour of $f(|Z_t|)$ until $Z_t$ reaches zero for the first time. Here on the other hand we will actually want to start both the marginal processes at the same point. Moreover, because of the modified construction of the coupling, which now behaves like a synchronous coupling for small values of $|Z_t|$, it can keep visiting zero infinitely often. A way to rigorously deal with this is to apply the version of the Meyer-It\^{o} formula that can be found e.g. as Theorem 71 in Chapter IV in \cite{protter}. We begin with computing the formula for $d|Z_t|$, by calculating $d|Z_t|^2$ first and then applying the It\^{o} formula once again to a~smooth approximation of the square root function, given e.g. by
\begin{equation*}
 S(r) := \begin{cases}
          -(1/8)\varepsilon^{-3/2} r^2 + (3/4) \varepsilon^{-1/2} r + (3/8) \varepsilon^{1/2} \,, &r < \varepsilon \,, \\
          \sqrt{r} \,, &r \geq \varepsilon \,.
         \end{cases}
\end{equation*}
A related argument was given by Zimmer in \cite{zimmer} in the context of infinite-dimensional diffusions, see Lemmas 2-5 therein. In our case, after two applications of the It\^{o} formula, we get
\begin{equation*}
 \begin{split}
  dS(|Z_t|^2) &= 4S'(|Z_t|^2)\lambda(Z_t)|Z_t|d\widetilde{W}_t + 2S'(|Z_t|^2) \langle Z_t, h_t + b(\widetilde{X}_t) - b(Y_t) \rangle dt \\
  &+ 4 S'(|Z_t|^2) \lambda^2(Z_t)dt + 8 S''(|Z_t|^2) \lambda^2(Z_t) |Z_t|^2 dt \\
  &+ \int_U \left( S\left(|Z_{t-} + g(\widetilde{X}_{t-},u) - g(Y_{t-},u)|^2\right) - S\left(|Z_{t-}|^2\right) \right) N(dt,du) \\
  &- 2 \int_U S'(|Z_{t-}|^2) \langle Z_{t-} , g(\widetilde{X}_{t-},u) - g(Y_{t-},u) \rangle \nu(du) dt \,.
 \end{split}
\end{equation*}
Since for any $r \in [0, \infty)$ we have $S(r) \to \sqrt{r}$ when $\varepsilon \to 0$, we can also show almost sure convergence of the integrals appearing in the formula above. For example, using the fact that $S$ is concave, for any $a$, $b \geq 0$ we have $S(a) - S(b) \leq S'(b)(b-a)$ and hence
\begin{equation*}
\begin{split}
 &\mathbb{E} \int_0^T \int_U \left( S\left(|Z_{t-} + g(\widetilde{X}_{t-},u) - g(Y_{t-},u)|^2\right) - S\left(|Z_{t-}|^2\right) \right) N(dt,du) \\
 &\leq \mathbb{E} \int_0^T \int_U S'(|Z_{t-}|^2) \left( |g(\widetilde{X}_{t-},u) - g(Y_{t-},u)|^2 + 2 \langle Z_{t-} , g(\widetilde{X}_{t-},u) - g(Y_{t-},u) \rangle \right) \nu(du) dt \,.
 \end{split}
\end{equation*}
Using the fact that $g$ is Lipschitz and that $\sup_{r \leq \varepsilon} S'(r) \lesssim \varepsilon^{-1/2}$, we see that the integral
\begin{equation*}
 \int_0^T \int_U \left( S\left(|Z_{t-} + g(\widetilde{X}_{t-},u) - g(Y_{t-},u)|^2\right) - S\left(|Z_{t-}|^2\right) \right) N(dt,du)
\end{equation*}
converges to
\begin{equation*}
 \int_0^T \int_U \mathbf{1}_{\{Z_{t-} \neq 0\}} \left( |Z_{t-} + g(\widetilde{X}_{t-},u) - g(Y_{t-},u)| - |Z_{t-}| \right) N(dt,du)
\end{equation*}
in $L^1$ and hence, via a subsequence, almost surely when $\varepsilon \to 0$. Dealing with the other integrals is even easier, cf. Lemmas 2 and 3 in \cite{zimmer} for analogous arguments. Thus we are able to get
\begin{equation*}
 \begin{split}
  d|Z_t| &= 2\lambda(Z_t) d\widetilde{W}_t + \mathbf{1}_{\{ Z_{t} \neq 0 \}} \frac{1}{|Z_{t}|} \langle Z_t, h_t + b(\widetilde{X}_t) - b(Y_t) \rangle dt \\
  &+ \int_U \mathbf{1}_{\{Z_{t-} \neq 0\}} \left( |Z_{t-} + g(\widetilde{X}_{t-},u) - g(Y_{t-},u)| - |Z_{t-}| \right) N(dt,du) \\
  &- \int_U \mathbf{1}_{\{Z_{t-} \neq 0\}} \frac{1}{|Z_{t-}|} \langle Z_{t-} , g(\widetilde{X}_{t-},u) - g(Y_{t-},u) \rangle \nu(du) dt \,.
 \end{split}
\end{equation*}
Now observe that the function $f$ defined by (\ref{defFEberle}) is twice continuously differentiable at all points except for $R_1$, whereas $f'$ exists and is continuous even at $R_1$. Therefore we can apply the Meyer-It\^{o} formula in its version given as Theorem 71 in Chapter IV in \cite{protter} to the process $(|Z_t|)_{t \geq 0}$ and the function $f$. For any $0 \leq s \leq r$ we get
\begin{equation}\label{ItoLocal}
 \begin{split}
  f(|Z_r|) - f(|Z_s|) &= 2 \int_s^r f'(|Z_{t}|) \lambda(Z_t) d\widetilde{W}_t \\
  &+ \int_s^r \mathbf{1}_{\{ Z_{t} \neq 0 \}} f'(|Z_{t}|)\frac{1}{|Z_{t}|} \langle Z_t, h_t + b(\widetilde{X}_t) - b(Y_t) \rangle dt \\
  &+ \int_s^r \int_U \mathbf{1}_{\{ Z_{t-} \neq 0 \}} f'(|Z_{t-}|) \frac{1}{|Z_{t-}|} \langle Z_{t-} , g(\widetilde{X}_{t-},u) - g(Y_{t-},u) \rangle \widetilde{N}(dt,du) \\
  &+ \int_s^r \int_U \mathbf{1}_{\{ Z_{t-} \neq 0 \}} \Bigg[ f(|Z_{t-} + g(X_{t-},u) - g(Y_{t-},u)|) - f(|Z_{t-}|) \\
 &- f'(|Z_{t-}|) \frac{1}{|Z_{t-}|} \langle Z_{t-} , g(X_{t-},u) - g(Y_{t-},u) \rangle \Bigg] N(dt,du) \\
&+ 2 \int_s^r f''(|Z_{t}|) \lambda^2(Z_t) dt \,.
 \end{split}
\end{equation}
We can see that the integrand in the integral with respect to $(\widetilde{W}_t)_{t \geq 0}$ in (\ref{ItoLocal}) is bounded (since $f'$ and $\lambda$ are bounded) and the integrand in the integral with respect to $\widetilde{N}$ is square integrable with respect to $\nu(du) dt$. Thus the expectations of both these integrals are zero. Moreover, the expectation of the integral with respect to $N$ in (\ref{ItoLocal}) can be dealt with in the same way as the expectation of the term $I_9$ in the proof of Theorem \ref{mainCouplingTheorem}, see (\ref{I9}). Thus, after taking the expectation everywhere in (\ref{ItoLocal}) and using the definition of $\kappa$, we get
\begin{equation}\label{ItoAfterExpectation}
 \begin{split}
  \mathbb{E}f(|Z_r|) - \mathbb{E}f(|Z_s|) &\leq \mathbb{E} \int_s^r \mathbf{1}_{\{ Z_{t} \neq 0 \}} | h_t | dt - \mathbb{E} \int_s^r \mathbf{1}_{\{ Z_{t} \neq 0 \}} f'(|Z_{t}|) |Z_{t}| \kappa(|Z_{t}|) dt\\
  &+ \mathbb{E} \int_s^r 2 f''(|Z_{t}|) \lambda^2(Z_t) dt \,.
 \end{split}
\end{equation}
We will now want to use the fact that the function $f$ defined by (\ref{defFEberle}) satisfies
\begin{equation}\label{differentialInequalityEberle}
 2f''(r) - r\kappa(r)f'(r) \leq -cf(r) \,.
\end{equation}
In particular, denoting $r_t := |Z_t|$, we get
\begin{equation*}
  2f''(r_t)\lambda^2(Z_t) - r_t \kappa(r_t)f'(r_t) \lambda^2(Z_t) + r_t \kappa(r_t)f'(r_t) - r_t \kappa(r_t)f'(r_t) \leq -cf(r_t)\lambda^2(Z_t) 
\end{equation*}
and thus
\begin{equation}\label{lemmaLambdaEstimates}
 - r_t \kappa(r_t)f'(r_t) + 2f''(r_t)\lambda^2(Z_t) \leq -cf(r_t)\lambda^2(Z_t)  + r_t \kappa(r_t)f'(r_t)(\lambda^2(Z_t) - 1) \,.
\end{equation}
Now observe that
\begin{equation*}
 -cf(r_t)\lambda^2(Z_t) = cf(r_t)(1 - \lambda^2(Z_t)) - cf(r_t) \leq c\delta - cf(r_t) \,, 
\end{equation*}
which holds since if $|Z_t| \geq \delta$, then $1 - \lambda^2(Z_t) = 0$ and if $|Z_t| \leq \delta$, then $1 - \lambda^2(Z_t) \leq 1$ and $cf(r_t) \leq c \delta$, which follow from the properties (\ref{lambdaProperties}) of the function $\lambda$ and the fact that $f(x) \leq x$ for any $x \in [0,\infty)$. Since obviously $-\kappa \leq \kappa^{-}$, we can further bound the right hand side of (\ref{lemmaLambdaEstimates}) by
\begin{equation}\label{lemmaLambdaEstimates2}
 c \delta - c f(r_t) + \kappa^{-}(r_t)(1 - \lambda^2(Z_t))r_t f'(r_t) \leq c \delta - c f(r_t) + \sup_{r \leq \delta} r \kappa^{-}(r) \,,
\end{equation}
where the last inequality follows from the fact that $1 - \lambda^2(Z_t) = 0$ when $|Z_t| \geq \delta$ and that $f' \leq 1$. If we denote
\begin{equation*}
 m(\delta) := c\delta + \sup_{r \leq \delta} r \kappa^{-}(r) \,,
\end{equation*}
then, from (\ref{lemmaLambdaEstimates}) and (\ref{lemmaLambdaEstimates2}) we obtain
\begin{equation}\label{lemmaLambdaEstimates3}
 - r_t \kappa(r_t)f'(r_t) + 2f''(r_t)\lambda^2(Z_t) \leq - c f(r_t) + m(\delta) \,.
\end{equation}
Hence, combining (\ref{ItoAfterExpectation}) with (\ref{lemmaLambdaEstimates3}) multiplied by $\mathbf{1}_{\{ r_t \neq 0\}}$, we get
\begin{equation*}
\begin{split}
  \mathbb{E}f(|Z_r|) - \mathbb{E}f(|Z_s|) &\leq -c \int_s^r \mathbb{E}\mathbf{1}_{\{ Z_{t} \neq 0 \}} f(|Z_t|) dt + \mathbb{E} \int_s^r \mathbf{1}_{\{ Z_{t} \neq 0 \}} \left( |h_t| + m(\delta) \right) dt\\
  &\leq -c \int_s^r \mathbb{E} f(|Z_t|) dt + \mathbb{E} \int_s^r |h_t| dt + \int_s^r m(\delta) dt \,.
\end{split}
\end{equation*}
Now observe that due to Assumption \ref{AssumptionD2}, we have $m(\delta) \to 0$ as $\delta \to 0$. We can also choose $\widetilde{X}_0 = Y_0$ so that $Z_0 = 0$. Eventually, applying the Gronwall inequality, we obtain
\begin{equation*}
  \mathbb{E}[f(|Z_t|)] \leq \mathbb{E} \int_0^t e^{c(s-t)} |h_s| ds \,,
\end{equation*}
which finishes the proof, since there exists a constant $C > 0$ such that for any $x \geq 0$ we have $x \leq Cf(x)$.
\end{proof}
\end{lemma}

\begin{proofDriftPerturbation}
 Once we have Lemma \ref{lemmaForMalliavin}, extending its result to the equation (\ref{twoNoisesSDE}) is quite straightforward. In comparison to the proof of Lemma \ref{lemmaForMalliavin}, the key step is to redefine $\kappa$ in order to include the additional coefficient $\sigma$ of the multiplicative Brownian noise (so that $\kappa$ satisfies (\ref{kappaSect4})) and then perform the same procedure as we did earlier (mixed reflection-synchronous coupling) only on the additive Brownian noise in order to construct processes like (\ref{XtildeNew}) and (\ref{MalliavinEstimateY}), where to the other noises we apply just the synchronous coupling. This way we can still get the inequality (\ref{differentialInequalityEberle}) with the same function $f$ as in the proof of Lemma \ref{lemmaForMalliavin}. The details are left to the reader, as they are just a repetition of what we have already presented. Once we obtain an inequality of the form (\ref{XYestimate}) for the equation (\ref{twoNoisesSDE}), we can use the Markov property of the process $(\widetilde{X}_t, Y_t)_{t \geq 0}$ to get (\ref{driftPerturbationInequality}), cf. Remark \ref{weakSolutionRemark}. 
 
 \qed
\end{proofDriftPerturbation}

\begin{proofGeneralMalliavin}
 
In order to keep notational simplicity, assume that we are dealing here with the equation $dX_t = b(X_t)dt + \sigma(X_t)dW_t$, i.e., the coefficient of the jump noise is zero. It does not influence our argument in any way, since we will need to perturb only the Gaussian noise. Recall that for a functional $f : \mathbb{R}^d \to \mathbb{R}$, the Malliavin derivative $\nabla_s f(X_t)$ is an $m$-dimensional vector $(\nabla_{s,1} f(X_t), \ldots, \nabla_{s,m} f(X_t))$, where $\nabla_{s,k} f(X_t)$ can be thought of as a derivative with respect to $(W^k_t)_{t \geq 0}$, where $W_t = (W^1_t, \ldots, W^m_t)$ is the driving $m$-dimensional Brownian motion.

We know that if $F$ is a random variable of the form $F = f(\int_0^T g^1_s dW_s, \ldots, \int_0^T g^N_s dW_s)$ for some smooth function $f : \mathbb{R}^N \to \mathbb{R}$ and $g^1, \ldots, g^N \in L^2([0,T]; \mathbb{R}^m)$ (i.e., $F \in \mathcal{S}$), then for any element $h \in H = L^2([0,T]; \mathbb{R}^m)$ we have
\begin{equation}\label{MalliavinDirectional}
 \langle \nabla F , h \rangle_{L^2([0,T]; \mathbb{R}^m)} = \int_0^t \langle \nabla_s F , h_s \rangle ds = \lim_{\varepsilon \to 0} \frac{1}{\varepsilon} \left( F(W_{\cdot} + \varepsilon \int_0^{\cdot} h_s ds) - F(W_{\cdot}) \right) \,,
\end{equation}
where convergence is in $L^2(\Omega)$. However, it is unclear whether (\ref{MalliavinDirectional}) holds also for arbitrary $F \in \mathbb{D}^{1,2}$ and in particular for $X_t$ (see the discussion in Appendix A in \cite{dinunno}, specifically Definitions A.10 and A.13). Nevertheless, for $F \in \mathbb{D}^{1,2}$ we can still prove that
\begin{equation}\label{MalliavinDirectionalExpectation}
 \mathbb{E} \langle \nabla F , h \rangle_{L^2([0,T]; \mathbb{R}^m)} = \lim_{\varepsilon \to 0} \frac{1}{\varepsilon} \mathbb{E} \left( F(W_{\cdot} + \varepsilon \int_0^{\cdot} h_s ds) - F(W_{\cdot}) \right) \,,
\end{equation}
even if we replace $h \in H$ with an adapted stochastic process $(\omega, t) \mapsto h_t(\omega)$ such that $\mathbb{E}\int_0^T |h_s|^2 ds < \infty$ and the Girsanov theorem applies (e.g. the Novikov condition for $h$ is satisfied).

Indeed, we know that for any $F \in \mathbb{D}^{1,2}$ and for any adapted square integrable $h$ we have
\begin{equation}\label{integrationByParts}
 \mathbb{E} \langle \nabla F , h \rangle_{L^2([0,T]; \mathbb{R}^m)} = \mathbb{E} \left[F \int_0^T h_s dW_s\right] \,.
\end{equation}
We recall now the proof of this fact, as we need to slightly modify it in order to get (\ref{MalliavinDirectionalExpectation}). As a reference, see e.g. Lemma A.15. in \cite{dinunno}, where (\ref{integrationByParts}) is proved only for $F \in \mathcal{S}$ and for deterministic $h$, but the argument can be easily generalized, or Theorems 1.1 and 1.2 in Chapter VIII of \cite{bass}. For now assume that $h$ is adapted and bounded (and thus it satisfies the assumptions of the Girsanov theorem). Then, starting from the right hand side of (\ref{integrationByParts}), we have
\begin{equation}\label{ibpProof}
 \begin{split}
  \mathbb{E} \left[F \int_0^T h_s dW_s\right] &= \mathbb{E}\left[F \frac{d}{d \varepsilon} \exp \left(\varepsilon \int_0^T h_s dW_s - \frac{1}{2}\varepsilon^2 \int_0^T |h_s|^2 ds \right)|_{\varepsilon = 0}\right] \\
  &= \mathbb{E} \left[F \lim_{\varepsilon \to 0} \frac{1}{\varepsilon} \left[\exp \left(\varepsilon \int_0^T h_s dW_s - \frac{1}{2}\varepsilon^2 \int_0^T |h_s|^2 ds \right) - 1\right]\right] \\
  &= \lim_{\varepsilon \to 0} \frac{1}{\varepsilon} \mathbb{E} \left[F \exp \left(\varepsilon \int_0^T h_s dW_s - \frac{1}{2}\varepsilon^2 \int_0^T |h_s|^2 ds \right) - F\right] \\
  &= \lim_{\varepsilon \to 0} \frac{1}{\varepsilon} \mathbb{E} \left[F(W_{\cdot} + \varepsilon \int_0^{\cdot} h_s ds) - F(W_{\cdot})\right] \,,
 \end{split}
\end{equation}
where in the last step we use the Girsanov theorem. In order to explain the third step, notice that the process
\begin{equation*}
 Z_t^{\varepsilon} := \exp \left(\varepsilon \int_0^t h_s dW_s - \frac{1}{2}\varepsilon^2 \int_0^t |h_s|^2 ds \right)
\end{equation*}
is the stochastic exponential of $\varepsilon \int_0^t h_s dW_s$ and thus it satisfies $dZ_t^{\varepsilon} = \varepsilon Z_s^{\varepsilon} h_s dW_s$, from which we get
\begin{equation*}
 \frac{1}{\varepsilon} \left[Z_T^{\varepsilon} - 1\right] = \int_0^T Z_s^{\varepsilon} h_s dW_s \,.
\end{equation*}
Now it is easy to see that since for any $\omega \in \Omega$ we have $Z_t^{\varepsilon}(\omega) \to 1$  with $\varepsilon \to 0$ and $Z_t^{\varepsilon}$ is uniformly bounded in $L^2(\Omega \times [0,T])$, there is a subsequence such that
\begin{equation*}
  \frac{1}{\varepsilon} \left[Z_T^{\varepsilon} - 1\right] = \int_0^T Z_s^{\varepsilon} h_s dW_s \to \int_0^T h_s dW_s \text{ as } \varepsilon \to 0 \,, \text{ in } L^2(\Omega) \,.
\end{equation*}
Thus the third step in (\ref{ibpProof}) holds for any $F \in L^2(\Omega)$ and in particular for any $F \in \mathbb{D}^{1,2}$. If $F$ is smooth, then the last expression in (\ref{ibpProof}) is equal to $\mathbb{E} \langle \nabla F , h \rangle_{L^2([0,T]; \mathbb{R}^m)}$, which proves (\ref{integrationByParts}) for any smooth $F$ and adapted, bounded $h$. Then (\ref{integrationByParts}) can be extended by approximation to any $F \in \mathbb{D}^{1,2}$ and any adapted, square integrable $h$. 

Now in order to prove (\ref{MalliavinDirectionalExpectation}), observe that the calculations in (\ref{ibpProof}) still hold when applied directly to an $F \in \mathbb{D}^{1,2}$ and an adapted, bounded $h$ (note that the argument does not work for general adapted, square integrable $h$ as we need to use the Girsanov theorem in the last step). Thus for any $F \in \mathbb{D}^{1,2}$ and any adapted, bounded $h$ we get
\begin{equation*}
 \mathbb{E} \langle \nabla F , h \rangle_{L^2([0,T]; \mathbb{R}^m)} = \mathbb{E} \left[F \int_0^T h_s dW_s\right] = \lim_{\varepsilon \to 0} \frac{1}{\varepsilon} \mathbb{E} \left[F(W_{\cdot} + \varepsilon \int_0^{\cdot} h_s ds) - F(W_{\cdot})\right] \,.
\end{equation*}
Since $X_t \in \mathbb{D}^{1,2}$ and $f$ is Lipschitz, we have $f(X_t) \in \mathbb{D}^{1,2}$ (cf. \cite{nualart}, Proposition 1.2.4), and hence
\begin{equation*}
 \mathbb{E} \langle \nabla f(X_t) , h \rangle_{L^2([0,T]; \mathbb{R}^m)} = \lim_{\varepsilon \to 0} \frac{1}{\varepsilon} \mathbb{E} \left( f(X_t)(W_{\cdot} + \varepsilon \int_0^{\cdot} h_s ds) - f(X_t)(W_{\cdot}) \right) 
\end{equation*}
holds for any adapted, bounded process $h$. From now on, we fix $t > 0$ and take $T = t$.

Recall that the process $(X_t)_{t \geq 0}$ is now given by $dX_t = b(X_t)dt + \sigma(X_t) dW_t$ and thus
\begin{equation}\label{perturbedDiffusion}
 X_t(W_{\cdot} + \varepsilon \int_0^{\cdot} h_s ds) = \int_0^t b(X_s)ds + \varepsilon \int_0^t \sigma(X_s) h_s ds + \int_0^t \sigma(X_s) dW_s \,.
\end{equation}
Hence, using the assumption (\ref{driftChangeBound}) from Theorem \ref{transportationInequalitiesTheorem} (taking $\varepsilon \sigma(X_s) h_s$ as the adapted change of drift and denoting the solution to (\ref{perturbedDiffusion}) by $(\widetilde{X}_t)_{t \geq 0}$) we obtain
\begin{equation*}
\begin{split}
 \mathbb{E} \left( f(X_t)(W_{\cdot} + \varepsilon \int_0^{\cdot} h_s ds) - f(X_t)(W_{\cdot}) \right) &= \mathbb{E} \left( f(X_t)(W_{\cdot} + \varepsilon \int_0^{\cdot} h_s ds) - f(Y'_t)(W_{\cdot}) \right) \\
 &\leq \varepsilon c_2(t) \mathbb{E} \int_0^t c_3(s) |\sigma(\widetilde{X}_s) h_s| ds \,, 
 \end{split}
\end{equation*}
where $\mathbb{E}f(X_t) = \mathbb{E}f(Y'_t)$, since $(X_t,Y'_t)_{t \geq 0}$ is a coupling. This in turn implies, together with our above calculations, that we have
\begin{equation}\label{MalliavinFirstBound}
 \mathbb{E} \langle \nabla f(X_t) , h \rangle_{L^2([0,t]; \mathbb{R}^m)} \leq c_2(t) \mathbb{E} \int_0^t  c_3(s) |\sigma(\widetilde{X}_s) h_s| ds \leq c_2(t) \sigma_{\infty} \mathbb{E} \int_0^t  c_3(s) |h_s| ds \,.
\end{equation}
Now by approximation we can show that the above inequality holds for any adapted process $h$ such that $\mathbb{E} \int_0^t |h_s|^2 ds < \infty$. Then, using the Cauchy-Schwarz inequality for $L^2(\Omega \times [0,t])$, we get
\begin{equation*}
 \mathbb{E} \langle \nabla f(X_t) , h \rangle_{L^2([0,t]; \mathbb{R}^m)} \leq   c_2(t) \sigma_{\infty} \left( \mathbb{E} \int_0^t c_3^2(s) ds \right)^{1/2} \left( \mathbb{E} \int_0^t |h_s|^2 ds \right)^{1/2} \,.
\end{equation*}
Moreover, observe that since $h$ is adapted, we have
\begin{equation}\label{EberleMalliavinTrick}
 \mathbb{E} \langle \nabla f(X_t) , h \rangle_{L^2([0,t]; \mathbb{R}^m)} = \mathbb{E} \int_0^t \langle \mathbb{E}[\nabla_s f(X_t) | \mathcal{F}_s] , h_s \rangle ds =: \mathbb{E} \langle  \mathbb{E}[\nabla_{\cdot} f(X_t) | \mathcal{F}_{\cdot}] , h_{\cdot} \rangle_{L^2([0,t]; \mathbb{R}^m)} \,.
\end{equation}
If we replace $h$ above with $h\sqrt{g}$ for some adapted, integrable, $\mathbb{R}_{+}$-valued process $g$, we get (by coming back to (\ref{MalliavinFirstBound}) and splitting $h$ and $\sqrt{g}$ via the Cauchy-Schwarz inequality)
\begin{equation*}
\begin{split}
 \mathbb{E} \langle \sqrt{g_{\cdot}} \mathbb{E}[\nabla_{\cdot} f(X_t) | \mathcal{F}_{\cdot}] , h_{\cdot} \rangle_{L^2([0,t]; \mathbb{R}^m)} &= \mathbb{E} \langle \mathbb{E}[\nabla_{\cdot} f(X_t) | \mathcal{F}_{\cdot}] , h_{\cdot}\sqrt{g_{\cdot}} \rangle_{L^2([0,t]; \mathbb{R}^m)}\\
 &\leq c_2(t) \sigma_{\infty} \left( \mathbb{E} \int_0^t g_s c_3^2(s) ds \right)^{1/2} \left( \mathbb{E} \int_0^t |h_s|^2 ds \right)^{1/2} \,.
 \end{split}
\end{equation*}
Since this holds for an arbitrary adapted, square integrable process $h$, we have
\begin{equation}\label{MalliavinBoundPreFinal}
 \mathbb{E} \int_0^t g_u |\mathbb{E}[ \nabla_u f(X_t) | \mathcal{F}_u]|^2 du \leq c^2_2(t)\sigma^2_{\infty} \mathbb{E} \int_0^t g_u c_3^2(u) du \,.
\end{equation}
Observe that in the inequality above we can integrate on any interval $[s,r] \subset [0,t]$. We can also approximate an arbitrary adapted, $\mathbb{R}_{+}$-valued process $g$ with processes $g \wedge n$ for $n \geq 1$, for which we have (\ref{MalliavinBoundPreFinal}). Then, by the Fatou lemma on the left hand side and the dominated convergence theorem on the right hand side, we get
\begin{equation*}
\begin{split}
 \mathbb{E} \int_s^r g_u |\mathbb{E}[ \nabla_u f(X_t) | \mathcal{F}_u]|^2 du &\leq \lim_{n \to \infty} \mathbb{E} \int_s^r (g_u \wedge n) |\mathbb{E}[ \nabla_u f(X_t) | \mathcal{F}_u]|^2 du \\
 &\leq c^2_2(t)\sigma^2_{\infty} \mathbb{E} \int_s^r g_u c_3^2(u) du \,.
 \end{split}
\end{equation*}
Hence we finally obtain (\ref{MalliavinExpectationEstimate}). In order to get (\ref{MalliavinLinftyEstimate}), we just need to go back to (\ref{MalliavinFirstBound}) and notice that it implies
\begin{equation*}
 \mathbb{E} \langle \nabla f(X_t) , h \rangle_{L^2([0,t]; \mathbb{R}^m)} \leq c_2(t) \sigma_{\infty} \sup_{u \leq t} c_3(u) \mathbb{E} \int_0^t |h_s| ds \,.
\end{equation*}
Since we can show that this holds for an arbitrary adapted $h$ from $L^1(\Omega \times [0,t])$, using (\ref{EberleMalliavinTrick}) and the fact that the dual of $L^1$ is $L^{\infty}$, we finish the proof.

\qed

\end{proofGeneralMalliavin}

\begin{proofCorollaryMalliavin}
 Note that from Theorem \ref{theoremDriftPerturbation} we obtain an inequality of the form (\ref{driftPerturbationInequality}), where on the right hand side we have either the coefficient $\sigma_1$ or $\sigma$, depending on whether we want to consider $\nabla^1$ or $\nabla^2$. Recall from the proof of Theorem \ref{theoremDriftPerturbation} that in order to get (\ref{driftPerturbationInequality}) we need to use the additive Brownian noise $(B_t^1)_{t \geq 0}$, regardless of which change of the drift we consider in the equation defining $(\widetilde{X}_t)_{t \geq 0}$ that appears therein. Therefore we need to assume $\det \sigma_1 > 0$ even if we are only interested in bounding the Malliavin derivative with respect to the multiplicative Brownian noise $(B_t^2)_{t \geq 0}$. Once we have (\ref{driftPerturbationInequality}), it is sufficient to apply Theorem \ref{generalMalliavinResult} with $c_2(t) = Ce^{-ct}$ and $c_3(s) = e^{cs}$. 
 
 \qed
 
\end{proofCorollaryMalliavin}

\subsection{Poissonian case}\label{MalliavinPoissonian}

Consider the solution $(X_t(x))_{t \geq 0}$ to
\begin{equation}\label{LevySDE}
 dX_t = b(X_t)dt + \sigma(X_t)dW_t + \int_U g(X_{t-},u) \widetilde{N}(dt,du)
\end{equation}
with initial condition $x \in \mathbb{R}^d$ as a functional of the underlying Poisson random measure $N = \sum_{j=1}^{\infty} \delta_{(\tau_j, \xi_j)}$. Then define
\begin{equation*}
 X^{(t,u)}(x) = X^{(t,u)}(x, N) := X(x, N + \delta_{(t,u)}) \,,
\end{equation*}
which means that we add a jump of size $g(X_{t-},u)$ at time $t$ to every path of $X$. Then
\begin{equation*}
 X^{(t,u)}_s(x) = X_s(x) \text{ for } s < t
\end{equation*}
and
\begin{equation*}
\begin{split}
 X^{(t,u)}_s(x) &= X_t(x) + g(X_{t-},u) + \int_t^s b(X^{(t,u)}_r(x)) dr \\
 &+ \int_t^s \sigma(X^{(t,u)}_r(x)) dW_r + \int_t^s \int_U g(X^{(t,u)}_{r-}(x), u) \widetilde{N}(dr,du) \text{ for } s \geq t \,.
\end{split}
\end{equation*}
This means that after time $t$, the process $(X^{(t,u)}_s(x))_{s \geq t}$ is a solution of the same SDE but with different initial condition, i.e., $X^{(t,u)}_t(x) = X_t(x) + g(X_{t-},u)$.

If the global dissipativity assumption is satisfied (like in \cite{lwu} and \cite{yma}), it is easy to show that the solution $(X_t)_{t \geq 0}$ to (\ref{LevySDE}) satisfies for any $x$ and $y \in \mathbb{R}^d$ the inequality
\begin{equation*}
 \mathbb{E} |X_t(x) - X_t(y)| \leq e^{-Kt}|x-y|
\end{equation*}
with some constant $K > 0$. Then we easily see that for any Lipschitz function $f : \mathbb{R}^d \to \mathbb{R}$ with $\| f \|_{\operatorname{Lip}} \leq 1$, if $t < T$ we have
\begin{equation*}
\begin{split}
 \mathbb{E} [ D_{t,u}f(X_T(x)) | \mathcal{F}_t] &\leq \mathbb{E} \left[ \left| f(X^{(t,u)}_T(x)) - f(X_T(x)) \right| | \mathcal{F}_t\right] \\
 &\leq \mathbb{E} \left[ \left| X^{(t,u)}_T(x) - X_T(x) \right| | \mathcal{F}_t \right] \\
 &\leq e^{-K(T-t)} |g(X_{t-},u)| \,.
\end{split}
\end{equation*}
In order to improve this result we will work under the assumption (\ref{initialChangeBound}) from Theorem \ref{transportationInequalitiesTheorem} stating that there exists a coupling $(X_t,Y_t)_{t \geq 0}$ of solutions to (\ref{LevySDE}) such that
\begin{equation}\label{section52couplingInequality}
 \mathbb{E}\left[\left|X_T - Y_T \right| / \mathcal{F}_t \right] \leq c_1(T-t)|X_t-Y_t|
\end{equation}
holds for any $T \geq t \geq 0$ with some function $c_1 : \mathbb{R}_{+} \to \mathbb{R}_{+}$. We fix $t > 0$ and we express the process $(X^{(t,u)}_s(x))_{s \geq 0}$ as
\begin{equation*}
 X^{(t,u)}_s(x) := \begin{cases}
                       X_s(x) & \text{ for } s < t \,, \\
                       \bar{X}_s & \text{ for } s \geq t \,,
                      \end{cases}
\end{equation*}
where $(\bar{X}_s)_{s \geq t}$ is a solution to (\ref{LevySDE}) started at $t$ with initial point $X_t(x) + g(X_{t-},u)$. Obviously both $(X_s)_{s \geq 0}$ and $(\bar{X}_s)_{s \geq t}$ have the same transition probabilities (since they are solutions to the same SDE satisfying sufficient conditions for uniqueness of its solutions in law). Thus we can apply our coupling to the process $(\bar{X}_s)_{s \geq t}$ to get a process $(\bar{Y}_s)_{s \geq t}$ with initial point $X_t(x)$ and the same transition probabilities as $(Y_s)_{s \geq 0}$ (and thus also $(X_s)_{s \geq 0}$). Now if we define the coupling time $\tau := \inf \{ r > t : \bar{X}_r = \bar{Y}_r \}$ then we can put
\begin{equation*}
  \widehat{Y}_s(x) := \begin{cases}
                       X_s(x) & \text{ for } s < t \,, \\
                       \bar{Y}_s & \text{ for } t \leq s < \tau \,, \\
                       \bar{X}_s & \text{ for }  s \geq \tau \,,
                      \end{cases} 
\end{equation*}
and we obtain a process with the same transition probabilities as $(X^{(t,u)}_s(x))_{s \geq 0}$ and thus also $(X_s(x))_{s \geq 0}$. This follows from a standard argument about gluing couplings at stopping times, see e.g. Subsection 2.2 in \cite{jwang4} for a possible approach. This way we get a~coupling $(X_s(x), \widehat{Y}_s(x) )_{s \geq 0}$ such that
\begin{equation*}
 \mathbb{E}\left[\left|X^{(t,u)}_T(x) - \widehat{Y}_T(x)\right| | \mathcal{F}_t \right] \leq c_1(T-t)|g(X_{t-},u)|
\end{equation*}
holds for any $T \geq t$ (from our construction we see that $X^{(t,u)}_t(x) - \widehat{Y}_t(x) = g(X_{t-},u)$ and we use (\ref{section52couplingInequality})). Now we can easily compute
\begin{equation}\label{MalliavinPoissonBound}
 \begin{split}
 \mathbb{E} [ D_{t,u}f(X_T(x)) | \mathcal{F}_t] &= \mathbb{E} [ f(X^{(t,u)}_T(x)) - f(X_T(x)) | \mathcal{F}_t]  \\
 &= \mathbb{E} [ f(X^{(t,u)}_T(x)) - f(\widehat{Y}_T(x)) | \mathcal{F}_t] \\
 &\leq \mathbb{E} \left[ \left| X^{(t,u)}_T(x) - \widehat{Y}_T(x) \right| | \mathcal{F}_t \right] \\
 &\leq c_1(T-t) |g(X_{t-},u)| \,,
\end{split}
\end{equation}
where we used the coupling property in the second step. In particular, if there exists a~measurable function $g_{\infty}: U \to \mathbb{R}$ such that $|g(x,u)| \leq g_{\infty}(u)$ for any $x \in \mathbb{R}^d$ and $u \in U$, then we obviously get
\begin{equation}\label{MalliavinPoissonBoundWithGinfty}
 \mathbb{E} [ D_{t,u}f(X_T(x)) | \mathcal{F}_t] \leq c_1(T-t) g_{\infty}(u) \,.
\end{equation}

To end this section, let us consider briefly the case of the equation (\ref{twoNoisesSDE}), where we have two jump noises, given by a L\'{e}vy process $(L_t)_{t \geq 0}$ and a Poisson random measure $N$. Then we can easily obtain analogous bounds on the Malliavin derivatives with respect to $(L_t)_{t \geq 0}$ and $N$, which we denote by $D^L$ and $D$, respectively. Namely, in the framework of Theorem \ref{mainCouplingTheorem} we obtain a coupling $(X_t,Y_t)_{t \geq 0}$ such that 
\begin{equation*}
 \mathbb{E}\left[\left|X_T - Y_T \right| / \mathcal{F}_t \right] \leq \widetilde{C} e^{-\widetilde{c}(T-t)}|X_t-Y_t|
\end{equation*}
holds for any $T \geq t \geq 0$ with some constants $\widetilde{C}$, $\widetilde{c} > 0$. Then, repeating the reasoning above, we easily get
\begin{equation}\label{twoNoisesMalliavinPoissonBound1}
  \mathbb{E} [ D^L_{t,u}f(X_T(x)) | \mathcal{F}_t] \leq \widetilde{C} e^{-\widetilde{c}(T-t)} u
\end{equation}
and
\begin{equation}\label{twoNoisesMalliavinPoissonBound2}
  \mathbb{E} [ D_{t,u}f(X_T(x)) | \mathcal{F}_t] \leq \widetilde{C} e^{-\widetilde{c}(T-t)} g_{\infty}(u) \,.
\end{equation}

\section{Proofs of transportation and concentration inequalities}\label{sectionProofsTransportation}

\begin{proofTransportationInequalities}
We first briefly recall the method of the proof of Theorem 2.2 in \cite{lwu} and its extension from \cite{yma} (however, we denote certain quantities differently from \cite{yma} to make the notation more consistent with the original one from \cite{lwu}). We will make use of the elements of Malliavin calculus described in Section 1. Specifically, we work on a probability space $(\Omega, \mathcal{F}, (\mathcal{F}_t)_{t \geq 0}, \mathbb{P})$ equipped with a Brownian motion $(W_t)_{t \geq 0}$ and a Poisson random measure $N$, on which we define the Malliavin derivative $\nabla$ with respect to $(W_t)_{t \geq 0}$ (a differential operator) and the Malliavin derivative $D$ with respect to $N$ (a difference operator). 
We use the Clark-Ocone formula, i.e., if $F$ is a functional such that the integrability condition (\ref{integrabilityConditionForClarkOcone}) is satisfied, then
\begin{equation}\label{ClarkOcone}
 F = \mathbb{E}F + \int_0^T \mathbb{E}[\nabla_t F | \mathcal{F}_t] dW_t + \int_0^T \int_U \mathbb{E}[D_{t,u} F |\mathcal{F}_t] \widetilde{N}(dt, du) \,.
\end{equation}
From the proof of Lemma 3.2 in \cite{yma} we know that if we show that there exists a deterministic function $h: [0,T] \times U \to \mathbb{R}$ such that $\int_0^T \int_U h(t,u)^2 \nu(du) dt < \infty$ and
\begin{equation}\label{MalliavinGeneralPoissonBound1}
 \mathbb{E}[D_{t,u} F |\mathcal{F}_t] \leq h(t,u)
\end{equation}
and there exists a deterministic function $j: [0,T] \to \mathbb{R}^m$ such that $\int_0^T |j(t)|^2 dt < \infty$ and
\begin{equation}\label{MalliavinGeneralBrownianBound1}
 |\mathbb{E}[\nabla_t F | \mathcal{F}_t]| \leq |j(t)| \,,
\end{equation}
then for any $C^2$ convex function $\phi: \mathbb{R} \to \mathbb{R}$ such that $\phi'$ is also convex, we have
\begin{equation}\label{convexInequality}
 \mathbb{E}\phi(F - \mathbb{E}F) \leq \mathbb{E} \phi \left( \int_0^T \int_U h(t,u) \widetilde{N}(dt, du) + \int_0^T j(t) dW_t \right) \,.
\end{equation}
In particular, for any $\lambda > 0$ we have
\begin{equation}\label{expInequality}
 \mathbb{E}e^{\lambda(F - \mathbb{E}F)} \leq \exp \left( \int_0^T \int_U (e^{\lambda h(t,u)} - \lambda h(t,u) - 1) \nu(du) dt + \int_0^T \frac{\lambda^2}{2} |j(t)|^2 dt \right) \,.
\end{equation}
The way to prove this is based on the forward-backward martingale method developed by Klein, Ma and Privault in \cite{klein}. On the product space $(\Omega^2, \mathcal{F}^2, \mathbb{P}^2)$ for any $(\omega, \omega') \in \Omega^2$ we can define
\begin{equation}\label{forwardMartingale}
 M_t(\omega, \omega') := \int_0^t \int_U \mathbb{E}[D_{s,u} F |\mathcal{F}_s](\omega) \widetilde{N}(\omega, ds, du) + \int_0^t \mathbb{E}[\nabla_s F | \mathcal{F}_s](\omega) dW_s(\omega) \,,
\end{equation}
which is a forward martingale with respect to the increasing filtration $\mathcal{F}_t \otimes \mathcal{F}$ on $\Omega^2$ and
\begin{equation}\label{backwardMartingale}
 M_t^{*}(\omega, \omega') := \int_t^T \int_U h(s,u) \widetilde{N}(\omega', ds, du) + \int_t^T j(s) dW_s(\omega') \,,
\end{equation}
which is a backward martingale with respect to the decreasing filtration $\mathcal{F} \otimes \mathcal{F}_t^{*}$, where $\mathcal{F}_t^{*}$ is the $\sigma$-field generated by $N([r, \infty),A)$ and $W_r$ for $r \geq t$ where $A$ are Borel subsets of $U$. Application of the forward-backward It\^{o} formula (see Section 8 in \cite{klein}) to $\phi(M_t + M_t^{*})$ and comparison of the characteristics of $M_t$ and $M_t^{*}$ shows that for any $s \leq t$ we have 
\begin{equation*}
 \mathbb{E}\phi(M_t + M_t^{*}) \leq \mathbb{E} \phi(M_s + M_s^{*}) \,.
\end{equation*}
This follows from Theorem 3.3 in \cite{klein}. However, it is important to note that if we replace (\ref{MalliavinGeneralBrownianBound1}) with a weaker assumption, stating that for any adapted, $\mathbb{R}_{+}$-valued process $g$ and for any $[s,r] \subset [0,T]$ we have
\begin{equation}\label{MalliavinGeneralBrownianBound2}
 \mathbb{E} \int_s^r g_u |\mathbb{E}[\nabla_u F | \mathcal{F}_u]|^2 du \leq \mathbb{E} \int_s^r g_u |j(u)|^2 du \,,
\end{equation}
then the argument from \cite{klein} still holds (check the page 493 in \cite{klein} and observe that what we need for the proof of Theorem 3.3 therein is that the integral of the process $\phi''(M_u + M_u^{*})$ appearing there is non-positive and that is indeed the case if $M$ and $M^{*}$ are given by (\ref{forwardMartingale}) and (\ref{backwardMartingale}), respectively, and the condition (\ref{MalliavinGeneralBrownianBound2}) holds).
Now we will use the fact that by the Clark-Ocone formula (\ref{ClarkOcone}) we know that $M_t + M_t^{*} \to F - \mathbb{E}F$ in $L^2$ as $t \to T$. Observe that since $\phi$ is convex, we have
\begin{equation*}
 \phi(M_t + M_t^{*}) - \phi(0) \geq \phi'(0)(M_t + M_t^{*}) 
\end{equation*}
and thus we can apply the Fatou lemma for $\phi(M_t + M_t^{*}) - \phi(0) - \phi'(0)(M_t + M_t^{*})$ to get
\begin{equation*}
 \mathbb{E}\phi(F - \mathbb{E}F) -\phi'(0)\mathbb{E}(M_T) \leq \lim_{t \to T} \mathbb{E}\phi(M_t + M_t^{*}) \,.
\end{equation*}
Here $\phi(0)$ cancels since it appears on both sides and by (\ref{ClarkOcone}) we know that $\mathbb{E}(M_T) = \mathbb{E}(F - \mathbb{E}F) = 0$. Thus we get
\begin{equation*}
 \mathbb{E}\phi(F - \mathbb{E}F) \leq \lim_{t \to T} \mathbb{E}\phi(M_t + M_t^{*}) \leq \lim_{t \to T} \mathbb{E}\phi(M_0^{*}) = \mathbb{E}\phi(M_0^{*}) \,,
\end{equation*}
which proves (\ref{convexInequality}).

Now we can return to the equation (\ref{generalSDEsect2}). Using the assumption (\ref{initialChangeBound}) we can get a~bound on the Malliavin derivative $D$ of a Lipschitz functional of $X_T(x)$, i.e., for any $f: \mathbb{R}^d \to \mathbb{R}$ with $\| f \|_{\operatorname{Lip}} \leq 1$ we have
\begin{equation}\label{MalliavinPoissonShortBound}
\mathbb{E} [ D_{t,u}f(X_T(x)) | \mathcal{F}_t] \leq c_1(T-t) |g(X_{t-},u)| \leq c_1(T-t) g_{\infty}(u)
\end{equation}
(see the discussion in Section \ref{MalliavinPoissonian}, in particular (\ref{MalliavinPoissonBound}) and (\ref{MalliavinPoissonBoundWithGinfty})). Note that the square integrability condition on the upper bound required in (\ref{MalliavinGeneralPoissonBound1}) is satisfied due to our assumptions on $g_{\infty}$. On the other hand, due to the assumption (\ref{driftChangeBound}), via Theorem \ref{generalMalliavinResult}, for any adapted $\mathbb{R}_{+}$-valued process $g$ and any $[s,r] \subset [0,T]$ we get
\begin{equation}\label{MalliavinFinalBound}
 \mathbb{E} \int_s^r g_u |\mathbb{E}[ \nabla_u f(X_T) | \mathcal{F}_u]|^2 du \leq c^2_2(T) \sigma_{\infty}^2 \mathbb{E} \int_s^r g_u c_3^2(u) du \,.
\end{equation}

It is easy to see that with our bounds, directly from (\ref{convexInequality}) we obtain (\ref{concentrationIneq1}). Note that as the integrand in the Brownian integral appearing in (\ref{concentrationIneq1}) we can take any $m$-dimensional function whose norm coincides with our upper bound in (\ref{MalliavinFinalBound}).
For the inequalities on the path space $\mathbb{D}([0,T];\mathbb{R}^d)$ we can still use our coupling $(X_s(x), \widehat{Y}_s(x) )_{s \geq 0}$ which we discussed in Section \ref{MalliavinPoissonian}. Denote by $\widehat{Y}_{[0,T]}$ a path of the process $(\widehat{Y}_s(x))_{t \in [0,T]}$. Then for any Lipschitz functional $F : \mathbb{D}([0,T];\mathbb{R}^d) \to \mathbb{R}$ (where we consider $\mathbb{D}([0,T];\mathbb{R}^d)$ equipped with the $L^1$ metric $d_{L^1}(\gamma_1, \gamma_2) := \int_0^T |\gamma_1(t) - \gamma_2(t)| dt$) such that $\| F \|_{\operatorname{Lip}} \leq 1$ we have
\begin{equation*}
  \begin{split}
   \mathbb{E}[D_{t,u}&F(X_{[0,T]}(x)) | \mathcal{F}_t] = \mathbb{E}[F(X^{(t,u)}_{[0,T]}(x)) - F(X_{[0,T]}(x))| \mathcal{F}_t] \\
   &= \mathbb{E}\left[ F(X^{(t,u)}_{[0,T]}(x)) - F(\widehat{Y}_{[0,T]}(x)) | \mathcal{F}_t \right] \leq \mathbb{E}\left[ \int_0^T  \left|X^{(t,u)}_r(x) - \widehat{Y}_r(x)\right| dr | \mathcal{F}_t \right] \\
   &=    \int_t^T \mathbb{E}\left[ \left|X^{(t,u)}_r(x) - \widehat{Y}_r(x)\right| | \mathcal{F}_t \right] dr \leq  \int_t^T c_1(r-t)|g(X_{t-},u)|  dr \\
   &\leq g_{\infty}(u) \int_t^T c_1(r-t) dr \,.
  \end{split}
\end{equation*}

In order to get a bound on $\mathbb{E}[\nabla_{\cdot} F(X_{[0,T]}(x)) | \mathcal{F}_{\cdot}]$, we proceed similarly as in the proof of Theorem \ref{generalMalliavinResult}, using again the coupling $(X_t, Y'_t)_{t \geq 0}$ satisfying the assumption (\ref{driftChangeBound}). Namely, we can show that for any bounded, adapted process $h$ we have
\begin{equation*}
 \begin{split}
  \mathbb{E} \langle \nabla F(&X_{[0,T]}(x)) , h \rangle_{L^2([0,T];\mathbb{R}^m)} \\
  &= \lim_{\varepsilon \to 0} \frac{1}{\varepsilon} \mathbb{E} \left( F(X_{[0,T]}(x))(W_{\cdot} + \varepsilon \int_0^{\cdot} h_u du ) - F(Y'_{[0,T]}(x))(W_{\cdot})\right) \\
  &\leq \lim_{\varepsilon \to 0} \frac{1}{\varepsilon} \int_0^T \mathbb{E} \left| X_r(x)(W_{\cdot} + \varepsilon \int_0^{\cdot} h_u du ) - Y'_r(x)(W_{\cdot}) \right| dr \\
  &\leq \int_0^T \left( c_2(r) \int_0^r c_3(u) \sigma_{\infty} |h_u| du \right) dr =  \int_0^T \left( \int_u^T c_2(r) c_3(u) \sigma_{\infty} |h_u| dr \right) du \\
  &\leq \left( \int_0^T \left( \int_u^T c_2(r) dr \right)^2 c_3^2(u) \sigma_{\infty}^2 du \right)^{1/2} \left( \int_0^T |h_u|^2 du \right)^{1/2} \,.
 \end{split}
\end{equation*}
Then we can extend this argument to obtain for any adapted $\mathbb{R}_{+}$-valued process $g$ and any $[s,t] \subset [0,T]$ the inequality
\begin{equation*}
  \mathbb{E} \int_s^t g_u |\mathbb{E}[ \nabla_u F(X_{[0,T]}(x)) | \mathcal{F}_u]|^2 du \leq \sigma_{\infty}^2 \mathbb{E} \int_s^t g_u c_3^2(u) \left( \int_u^T c_2(r) dr \right)^2 du \,.
\end{equation*}

This, due to (\ref{convexInequality}), gives (\ref{concentrationIneq2}). This finishes the proof of Theorem \ref{concentrationInequalitiesTheorem}. Notice that the inequalities therein are true even if the expectation on the right hand side is infinite. However, if we want to obtain transportation inequalities from Theorem \ref{transportationInequalitiesTheorem}, we need the Assumption \ref{AssumptionE}. Then we can apply our reasoning and the inequality (\ref{convexInequality}) with the function $\phi(x) = \exp (\lambda x)$ and after simple calculations we obtain (\ref{expInequality}), which in the case of our bounds on Malliavin derivatives reads as
\begin{equation}\label{preTransIneq1}
  \mathbb{E}e^{\lambda (f(X_T(x)) - p_Tf(x))} \leq \exp \left( \int_0^T \beta(\lambda c_1(T-t)) dt + \frac{\lambda^2}{2} \sigma_{\infty}^2 c_2^2(T) \int_0^T c_3^2(t) dt \right)
\end{equation}
and on the path space as
\begin{equation}\label{preTransIneq2}
\begin{split}
 \mathbb{E}&e^{\lambda \left( F(X_{[0,T]}(x)) - \mathbb{E}F(X_{[0,T]}(x)) \right)} \\
 &\leq \exp \left( \int_0^T \beta\left(\lambda \int_t^T c_1(r-t)dr \right) dt + \frac{\lambda^2}{2} \sigma_{\infty}^2 \int_0^T c_3^2(t) \left(\int_t^T c_2(r)dr \right)^2 dt \right) \,.
 \end{split}
\end{equation}
Then, by the Gozlan-L\'{e}onard characterization (\ref{gozlanCharacterization}) and the Fenchel-Moreau theorem, we easily get (\ref{transIneq1}) from (\ref{preTransIneq1}) and (\ref{transIneq2}) from (\ref{preTransIneq2}).

\qed

\end{proofTransportationInequalities}

\begin{remark}\label{finiteIntensityRemark}
Note that if instead of (\ref{initialChangeBound}) we have an inequality like
\begin{equation}\label{alternativeCouplingInequality}
\mathbb{E}[|X_t - Y_t| / \mathcal{F}_s ]\leq c_1(t-s) (|X_s-Y_s| + 1) \,,
\end{equation}
then, by the same reasoning as in Section \ref{MalliavinPoissonian}, instead of (\ref{MalliavinPoissonShortBound}) we get
\begin{equation*}
\mathbb{E} [ D_{t,u}f(X_T(x)) | \mathcal{F}_t] \leq c_1(T-t) (g_{\infty}(u) + 1) \,.
\end{equation*}
Then, if we want to obtain transportation or concentration inequalities, $ g_{\infty}(u) + 1$ has to be square integrable with respect to the measure $\nu$. However, if $\nu$ is a L\'{e}vy measure, this implies that $\nu$ has to be finite. This could still allow us to obtain some interesting results in certain cases that are not covered by Corollary \ref{corollaryTransportation}, where Assumption \ref{Assumption6} is required, which we do not need to obtain (\ref{alternativeCouplingInequality}) (cf. Remark \ref{finiteLevyMeasureRemark}). For the sake of brevity, we skip the details.
\end{remark}

\begin{proofCorollaryTransportation}
 In the presence of two Gaussian and two jump noises, we use the Clark-Ocone formula of the form
 \begin{equation*}
 \begin{split}
  F &= \mathbb{E}F + \int_0^T \mathbb{E}[\nabla^1_t F | \mathcal{F}_t] dB^1_t + \int_0^T \mathbb{E}[\nabla^2_t F | \mathcal{F}_t] dB^2_t \\
  &+ \int_0^T \int_U \mathbb{E}[D^L_{t,u} F |\mathcal{F}_t] \widetilde{N^L}(dt, du) + \int_0^T \int_U \mathbb{E}[D_{t,u} F |\mathcal{F}_t] \widetilde{N}(dt, du) \,,
  \end{split}
 \end{equation*}
which holds for square integrable functionals $F$, where $\nabla^1$, $\nabla^2$, $D^L$ and $D$ are the Malliavin derivatives with respect to $(B^1_t)_{t \geq 0}$, $(B^2_t)_{t \geq 0}$, $N^L$ and $N$, respectively (see e.g. Theorem 12.20 in \cite{dinunno}). Then we proceed as in the proof of Theorem \ref{transportationInequalitiesTheorem}, using the fact that under our assumptions, Theorem \ref{mainCouplingTheorem} and Theorem \ref{theoremDriftPerturbation} provide us with couplings such that the conditions (\ref{initialPerturbationInequality}) and (\ref{driftPerturbationInequality}) are satisfied and this allows us to obtain the required bounds on the Malliavin derivatives (of the type (\ref{MalliavinPoissonShortBound}) and (\ref{MalliavinFinalBound})). More precisely, under our assumptions we obtain (\ref{additiveMalliavinEstimate}) and (\ref{multiplicativeMalliavinEstimate}) from Corollary \ref{corollaryMalliavinBounds}, whereas (\ref{twoNoisesMalliavinPoissonBound1}) and (\ref{twoNoisesMalliavinPoissonBound2}) follow from our reasoning at the end of Section \ref{MalliavinPoissonian}. Combining all these bounds and using (\ref{convexInequality}), just like in the proof of Theorem \ref{transportationInequalitiesTheorem}, allows us to obtain the desired transportation inequalities. Furthermore, taking $T \to \infty$ in the $\alpha_T$-$W_1H$ inequality, we obtain (\ref{alphaW1HforInvariant}) by the argument from the proof of Lemma 2.2 in \cite{djellout}.
 
 \qed
\end{proofCorollaryTransportation}

\section*{Acknowledgement}
I would like to thank Arnaud Guillin and Liming Wu for suggesting the topic and for their hospitality and fruitful discussions during my visit to Clermont-Ferrand. This visit was financed by DAAD and took place while I was a PhD student funded by the Bonn International Graduate School of Mathematics. I am also grateful to my PhD advisor, Andreas Eberle, for many helpful suggestions and constant support. This research was partially supported by the ERC grant no. 694405.


\begin{thebibliography}{99}

\bibitem[1]{albeverio} S. Albeverio, Z. Brze\'{z}niak, J.-L. Wu, Existence of global solutions and invariant measures for stochastic differential equations driven by Poisson type noise with non-Lipschitz coefficients, J. Math. Anal. Appl. 371 (2010), no. 1, 309-322.

\bibitem[2]{apple} D. Applebaum, L\'{e}vy Processes and Stochastic Calculus, 2nd ed., Cambridge University Press, 2009.

\bibitem[3]{bakryemery} D. Bakry, M. \'{E}mery, Diffusions hypercontractives, S\'{e}minaire de probabilit\'{e}s, XIX, 1983/84, 177-206,
Lecture Notes in Math., 1123, Springer, Berlin, 1985. 

\bibitem[4]{basscranston} R. F. Bass, M. Cranston, The Malliavin calculus for pure jump processes and applications to local time, Ann. Probab. 14 (1986), no. 2, 490-532.

\bibitem[5]{bass} R. F. Bass, Diffusions and Elliptic Operators, Springer-Verlag, New York, 1998.

\bibitem[6]{bichteler} K. Bichteler, J.-B. Gravereaux, J. Jacod, Malliavin calculus for processes with jumps, Stochastics Monographs, 2. Gordon and Breach Science Publishers, New York, 1987.

\bibitem[7]{bismut} J.-M. Bismut, Calcul des variations stochastique et processus de sauts, Z. Wahrsch. Verw. Gebiete 63 (1983), no. 2, 147-235. 

\bibitem[8]{bobkov} S. G. Bobkov, F. G\"{o}tze, Exponential integrability and transportation cost related to logarithmic Sobolev inequalities, J. Funct. Anal. 163 (1999), no. 1, 1-28. 

\bibitem[9]{carmona} R. Carmona, W. C. Masters, B. Simon, Relativistic Schr\"{o}dinger operators: asymptotic behavior of the eigenfunctions, J. Funct. Anal. 91 (1990), no. 1, 117-142. 

\bibitem[10]{dinunno} G. Di Nunno, B. {\O}ksendal, F. Proske, Malliavin Calculus for L\'{e}vy Processes with Applications to Finance, Universitext. Springer-Verlag, Berlin, 2009.

\bibitem[11]{djellout} H. Djellout, A. Guillin, L. Wu, Transportation cost-information inequalities and applications to random dynamical systems and diffusions, Ann. Probab. 32 (2004), no. 3B, 2702-2732. 

\bibitem[12]{eberle} A. Eberle, Reflection couplings and contraction rates for diffusions, Probab. Theory Related Fields 166 (2016), no. 3-4, 851-886.

\bibitem[13]{eberleguillin} A. Eberle, A. Guillin, R. Zimmer, Quantitative Harris type theorems for diffusions and McKean-Vlasov processes, Trans. Amer. Math. Soc. (2018), in press, \url{https://doi.org/10.1090/tran/7576}.

\bibitem[14]{gozlan} N. Gozlan, C. L\'{e}onard, A large deviation approach to some transportation cost inequalities, Probab. Theory Related Fields 139 (2007), no. 1-2, 235-283.

\bibitem[15]{gozlansurvey} N. Gozlan, C. L\'{e}onard, Transport inequalities. A survey, Markov Process. Related Fields 16 (2010), no. 4, 635-736. 

\bibitem[16]{guillinleonard} A. Guillin, C. L\'{e}onard, L. Wu, N. Yao, Transportation-information inequalities for Markov processes, Probab. Theory Related Fields 144 (2009), no. 3-4, 669-695.

\bibitem[17]{gyongykrylov} I. Gy\"{o}ngy, N. V. Krylov, On stochastic equations with respect to semimartingales. I., Stochastics 4 (1980/81), no. 1, 1-21. 

\bibitem[18]{huang} L. Huang, Density estimates for SDEs driven by tempered stable processes, preprint, arXiv:1504.04183.

\bibitem[19]{klein} T. Klein, Y. Ma, N. Privault, Convex concentration inequalities and forward-backward stochastic calculus, Electron. J. Probab. 11 (2006), no. 20, 486-512.

\bibitem[20]{komorowski} T. Komorowski, A. Walczuk, Central limit theorem for Markov processes with spectral gap in the Wasserstein metric, Stochastic Process. Appl. 122 (2012), no. 5, 2155-2184.

\bibitem[21]{lastpenrose} G. Last, M. Penrose, Martingale representation for Poisson processes with applications to minimal variance hedging, Stochastic Process. Appl. 121 (2011), no. 7, 1588-1606.

\bibitem[22]{lindvall} T. Lindvall, L. C. G. Rogers, Coupling of multidimensional diffusions by reflection, Ann. Probab. 14 (1986), no. 3, 860-872.

\bibitem[23]{luowang} D. Luo, J. Wang, Refined basic couplings and Wasserstein-type distances for SDEs with L\'{e}vy noises, Stochastic Process. Appl. (2018), in press, \url{https://doi.org/10.1016/j.spa.2018.09.003}.

\bibitem[24]{lokka} A. L{\o}kka, Martingale representation of functionals of L\'{e}vy processes, Stochastic Anal. Appl. 22 (2004), no. 4, 867-892.

\bibitem[25]{yma} Y. Ma, Transportation inequalities for stochastic differential equations with jumps, Stochastic Process. Appl. 120 (2010), no. 1, 2-21.

\bibitem[26]{maprivault} Y. Ma, N. Privault, Convex concentration for some additive functionals of jump stochastic differential equations,  Acta Math. Sin. (Engl. Ser.) 29 (2013), no. 8, 1449-1458. 

\bibitem[27]{majka} M. B. Majka, Coupling and exponential ergodicity for stochastic differential equations driven by L\'{e}vy processes, Stochastic Process. Appl. 127 (2017), no. 12, 4083-4125.

\bibitem[28]{nualart} D. Nualart, The Malliavin Calculus and Related Topics, second ed., Springer-Verlag, Berlin, 2006.

\bibitem[29]{ottovillani} F. Otto, C. Villani, Generalization of an inequality by Talagrand and links with the logarithmic Sobolev inequality, J. Funct. Anal. 173 (2000), no. 2, 361-400. 

\bibitem[30]{picardfrench} J. Picard, Formules de dualit\'{e} sur l'espace de Poisson, Ann. Inst. H. Poincar\'{e} Probab. Statist. 32 (1996), no. 4, 509-548.

\bibitem[31]{picard} J. Picard, On the existence of smooth densities for jump processes, Probab. Theory Related Fields 105 (1996), no. 4, 481-511. 

\bibitem[32]{priola} E. Priola, F. Y. Wang, Gradient estimates for diffusion semigroups with singular coefficients, J. Funct. Anal. 236 (2006), no. 1, 244-264. 

\bibitem[33]{protter} P. Protter, Stochastic Integration and Differential Equations, Second edition. Version 2.1., Springer-Verlag, Berlin, 2005.

\bibitem[34]{sato} K. Sato, L\'{e}vy Processes and Infinitely Divisible Distributions, Cambridge University Press, 1999.

\bibitem[35]{ryznar} M. Ryznar, Estimates of Green function for relativistic $\alpha$-stable process, Potential Anal. 17 (2002), no. 1, 1-23. 

\bibitem[36]{shao} J. Shao, C. Yuan, Transportation-cost inequalities for diffusions with jumps and its application to regime-switching processes, J. Math. Anal. Appl. 425 (2015), no. 2, 632-654. 

\bibitem[37]{villani} C. Villani, Optimal Transport. Old and New., Springer-Verlag, Berlin, 2009. 

\bibitem[38]{jwang4} J. Wang, $L^p$-Wasserstein distance for stochastic differential equations driven by L\'{e}vy processes, Bernoulli 22 (2016), no. 3, 1598-1616.

\bibitem[39]{lwu} L. Wu, Transportation inequalities for stochastic differential equations of pure jumps, Ann. Inst. Henri Poincar\'{e} Probab. Stat. 46 (2010), no. 2, 465-479.

\bibitem[40]{zimmer} R. Zimmer, Explicit contraction rates for a class of degenerate and infinite-dimensional diffusions, Stoch. Partial Differ. Equ. Anal. Comput. 5 (2017), no. 3, 368-399.



\end{thebibliography}
\end{document}